\documentclass[12pt,reqno]{amsart}
\usepackage{amsfonts,amsthm}

 \usepackage{hyperref}

\input epsf.sty

\newtheorem{theorem}{Theorem}[section]
\newtheorem{lemma}[theorem]{Lemma}
\newtheorem{proposition}[theorem]{Proposition}
\newtheorem{corollary}[theorem]{Corollary}

\newtheorem{definition}{Definition}[section]

\theoremstyle{remark}

\numberwithin{equation}{section}

\newcommand{\eps}{\epsilon}
\newcommand{\cal}{\mathcal}
\newcommand{\half}{\frac{1}{2}}
\renewcommand{\S}{\cal S}
\renewcommand{\b}{\beta}
\newcommand{\db}{\beta}

\newcommand{\vlim}{\tilde{u_\eps}}
\newcommand{\cefrac}{\frac{1+\theta(\eps)}{2}}
\newcommand{\fracce}{\frac{1-\theta(\eps)}{2}}
\newcommand{\epsfrac}{\frac{1+\theta(\eps)}{1-\theta(\eps)}}

\newcommand{\R}{\mathbb{R}}
\newcommand{\Z}{\mathbb{Z}}
\newcommand{\N}{\mathbb{N}}

\newcommand{\p}{\partial}
\def\sgn{\mathrm{sgn}}
\def\1{\mathbf{1}}

\newcommand{\be}{\begin{eqnarray}}
\newcommand{\ee}{\end{eqnarray}}

\begin{document}

%
%

\title[Biased tug-of-war and the biased infinity Laplacian]{Biased tug-of-war, the biased infinity Laplacian,\\ and comparison with exponential cones}

\author[Y. Peres]{Yuval Peres}
\address{Microsoft Research, One Microsoft Way, Redmond, WA 98052-6399, \and Department of Statistics, University of California, Berkeley, CA 94720}
\email{peres@microsoft.com}

\author[G. Pete]{G\'abor Pete}
\address{Department of Mathematics, University of Toronto, 40 St George St., Toronto, ON, M5S 2E4, Canada}
\email{gabor@math.toronto.edu}

\author[S. Somersille]{Stephanie Somersille}
\address{Department of Mathematics, University of California, Berkeley, CA 94720-3840}
\email{steph@math.berkeley.edu}





\begin{abstract}
We prove that if $U\subset \R^n$ is an open domain whose closure $\overline U$ is compact in the path metric, and $F$ is a Lipschitz function on $\partial U$, then for each $\b \in \R$ there exists a unique viscosity solution to the $\b$-biased infinity Laplacian equation
$$\b |\nabla u| + \Delta_\infty u=0$$
on $U$ that extends $F$, where $\Delta_\infty u= |\nabla u|^{-2} \sum_{i,j} u_{x_i}u_{x_ix_j} u_{x_j}$.

In the proof, we extend the tug-of-war ideas of Peres, Schramm, Sheffield and Wilson, and define the $\beta$-biased $\eps$-game as follows. The starting position is $x_0 \in U$. At the $k^\text{th}$ step the two players toss a suitably biased coin (in our key example, player I wins with odds of $\exp(\b\eps)$ to 1), and the winner chooses $x_k$ with $d(x_k,x_{k-1}) < \eps$. The game ends when $x_k \in \partial U$, and player II pays the amount $F(x_k)$ to player I. We prove that the value $u^{\eps}(x_0)$ of this game exists, and that $\|u^\eps - u\|_\infty \to 0$ as $\eps \to 0$, where $u$ is the unique extension of $F$ to $\overline{U}$ that satisfies comparison with $\db$-exponential cones. Comparison with exponential cones is a notion that we introduce here, and generalizing a theorem of Crandall, Evans and Gariepy regarding comparison with linear cones, we show that a continuous function satisfies comparison with $\db$-exponential cones if and only if it is a viscosity solution to the $\b$-biased infinity Laplacian equation.
\end{abstract}



\maketitle

%
%
%
%
\section{Introduction}\label{s.intro}

Informally, the {\bf infinity Laplacian equation} $\Delta_\infty u=0$ is the properly interpreted Euler-Lagrange equation associated with minimizing the functional $(u,X) \mapsto \|\nabla u\|_{L^\infty(X)}$ for $X\subseteq\R^n$, and the operator
\be\label{inftyop}
\Delta_\infty u := |\nabla u|^{-2} \sum_{i,j=1}^n u_{x_i}u_{x_ix_j} u_{x_j}
\ee
is the second derivative in the direction of the gradient $\nabla u=(u_{x_i})_{i=1}^n$. It was introduced and studied by Aronsson in 1967 \cite{A1,A2}, but even the basic existence and uniqueness questions have proven difficult, largely because of the non-smoothness of solutions. Several approaches were developed to overcome this problem, including the notion of viscosity solutions (see \cite{CIL}) and the method of comparison with cones, developed by Crandall, Evans and Gariepy \cite{CEG}. It was only in 1993 that Jensen \cite{Jensen} showed that if $U\subset \R^n$ is a bounded domain, and $F$ is a Lipschitz function on $\p U$, then a continuous function $u$ is a viscosity solution to $\Delta_\infty u=0$ if and only if it is a so-called absolutely minimizing Lipschitz extension of $F$. Jensen also proved uniqueness in this setting. See \cite{ACJ} for a survey of absolutely minimizing Lipschitz extensions. Then, in their recent paper \cite{PSSW2}, Peres, Schramm, Sheffield and Wilson introduced a new perspective by applying game theory to these problems. Using the game random-tug-of-war, they proved the most general existence and uniqueness results to date for solving equations involving the operator $\Delta_\infty$. See Section~\ref{s.open} for more references.

In the present work, we consider an extension of the infinity Laplacian PDE, the {\bf biased version}
\be\label{biaslap}
\Delta_\infty u + \b |\nabla u| =0,\qquad \b\in\R\text{ fixed,}
\ee
which is very natural both from game theory and PDE points of view. It appears, e.g., in \cite{BEJ}, which considers generalizations of~(\ref{inftyop}) that are motivated by game theory, but does not prove any existence or uniqueness results. On the other hand,~(\ref{biaslap}) is covered by a general existence and uniqueness theorem of Barles and Busca \cite{BB}, proved by PDE methods under strong smoothness conditions on the boundary of the domain. Our work gives a much better understanding of the corresponding game (the {\bf $\b$-biased $\eps$-tug-of-war}, where the coin is suitably biased) than \cite{BEJ}, and, using this, proves a more general existence and uniqueness result for~(\ref{biaslap}) than \cite{BB}. Our paper follows the general strategy of \cite{PSSW2}, but there are several new elements and difficulties that were not present in the unbiased case, partly because the connection to absolutely minimizing Lipschitz extensions is not there any more. The link between biased tug-of-war and the PDE is now given by a new characterization of viscosity solutions of~(\ref{biaslap}) in terms of {\bf comparison with exponential cones}. This comparison property also makes it possible to state and prove extension results beyond the setting of $\R^n$, namely, in general length spaces. We now describe our results in detail.
\medskip

For any metric space $X$, we can define a new metric, the {\bf path metric} $d(x_1,x_2)$, as the infimum of the lengths of all continuous paths from $x_1$ to $x_2$ which are contained in $X$. A metric space is called a {\bf length space} if it is equipped with this path metric. Note that the closure of a bounded domain $U \subset \mathbb{R}^n$ is not necessarily compact in the path metric; an example is an infinite inward spiral. On the other hand, the union of a finite number of convex closed sets in $\R^n$ is a compact length space.

Given a length space $X$, some $Y\subset X$, a real function $F$ on $Y$, a small $\eps>0$, and any $\rho(\eps)>0$, the {\bf biased $\eps$-tug-of-war} game is defined as follows.
The starting position is $x_0 \in X\backslash Y$. At the $k^\text{th}$ step the two players toss a biased coin which player I wins with odds of $\rho(\eps)$ to 1, i.e., with probability $\rho(\eps)/(\rho(\eps)+1)$, and the winner chooses $x_k$ with $d(x_k,x_{k-1}) < \eps$. (Even if $X\subset \R^n$, the moves are measured in the path metric of $X$, not in the Euclidean metric.) The game ends when $x_k \in Y$, and player II pays the amount $F(x_k)$ to player I. Assuming that the value of this game exists (defined and proved below), it is denoted by $u^{\eps}(x_0)$.

In particular, for any fixed $\b \in \R$, a {\bf $\b$-biased $\eps$-tug-of-war} game is defined by assuming that the odds function satisfies the following:
\begin{equation}
\label{e.rhocond}
\begin{aligned}
\rho(x)&\text{ is a smooth function with $\rho(0)=1$ and $\rho'(0)=\beta$,}\\
&\text{ whose logarithm is concave or convex.}
\end{aligned}
\end{equation}
Note that if the probability for player I to win a coin toss is $\frac{1+\theta(\eps)}{2}$, then we have
\begin{equation}
\rho(\eps)=\frac{1+\theta(\eps)}{1-\theta(\eps)}\qquad\hbox{and}\qquad \theta(\eps)=\frac{\rho(\eps)-1}{\rho(\eps)+1}.
\label{e.rhotheta}
\end{equation}
A natural example is $\theta(\eps)=\b\eps/2$,
but for reasons that will be explained shortly, the most important example for us will be $$\rho_0(\eps):=e^{\b\eps}\qquad\hbox{and}\qquad
\theta_0(\eps):=\tanh\frac{\b\eps}{2}=\frac{e^{\b\eps}-1}{e^{\b\eps}+1}=\frac{\b\eps}{2} - \frac{(\b\eps)^3}{24} + O(\b\eps)^5\,.
$$

Our main result is the following. For the definition of viscosity solutions, see Section~\ref{s.equiv}.

\begin{theorem}
\label{t.main}
Let $U\subset \R^n$ be an open domain with closure $X=\overline{U}$ compact in the path metric, and let $F$ be a Lipschitz function on $Y=\partial U$. Fix $\b \in \R$. Then there exists a unique extension of $F$ to $X$ that is a viscosity solution of~(\ref{biaslap}) on $U$. This solution is the uniform limit of the value functions $u^\eps(x)$ for the $\b$-biased $\eps$-tug-of-war games on $X$.
\end{theorem}

Let us point out that, unlike for the unbiased equation, the compactness assumption on $X$ cannot be dropped. For instance, equation~(\ref{biaslap}) in $\R$ takes the form $u''(x)+\db|u'(x)|=0$, so, if $X=[0,\infty)$, $Y=\{0\}$ and $F(0)=0$, then the functions $u(x)=\sgn(\b)A(1-e^{-|\db| x})$ with $A\in\R_{\geq 0}$ are all bounded extensions on $X$ solving~(\ref{biaslap}). These one-dimensional $\db$-exponential functions will in fact play a crucial role in the sequel.

\begin{definition} \label{defEC}
Given $\b\in\R$, $A\geq 0$, and $x_0\in X$ in a length space, we define
\begin{eqnarray}
\label{expcone}
C^+_{x_0}(x) &:=&\sgn(\b)\,A\left(1-e^{-\b\, d(x_0,x)}\right) \geq 0,\\
C^-_{x_0}(x) &:=&\sgn(\b)\,A\left(1-e^{\b\,d(x_0,x)}\right) \leq 0\,.\nonumber
\end{eqnarray}
For any $B\in \R$, the functions $C^+_{x_0}(x) + B$ and $C^-_{x_0}(x) +B$  are called {positive} and {negative} {$\b$-exponential cones} centered at $x_0$ with ``slope'' $A$.
\end{definition}

\begin{figure}[htbp]
\centerline{\raise 2.5 cm \hbox{\epsfxsize=0.4\textwidth \epsffile{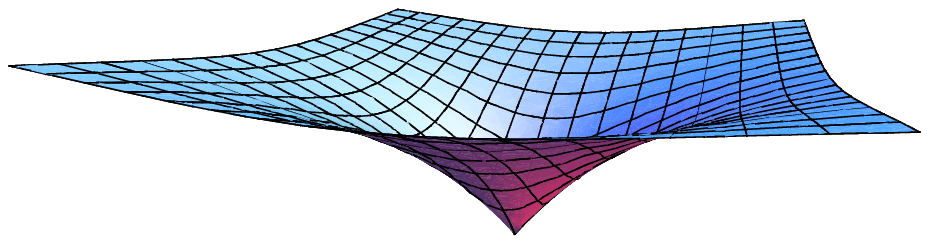}}
\hskip 0.5 cm
\raise 0.5 cm \hbox{\epsfxsize=0.17\textwidth \epsffile{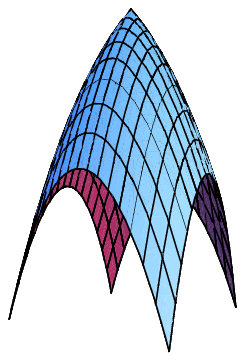}}
\hskip 0.5 cm
\raise 0.3 cm \hbox{\epsfxsize=0.30\textwidth \epsffile{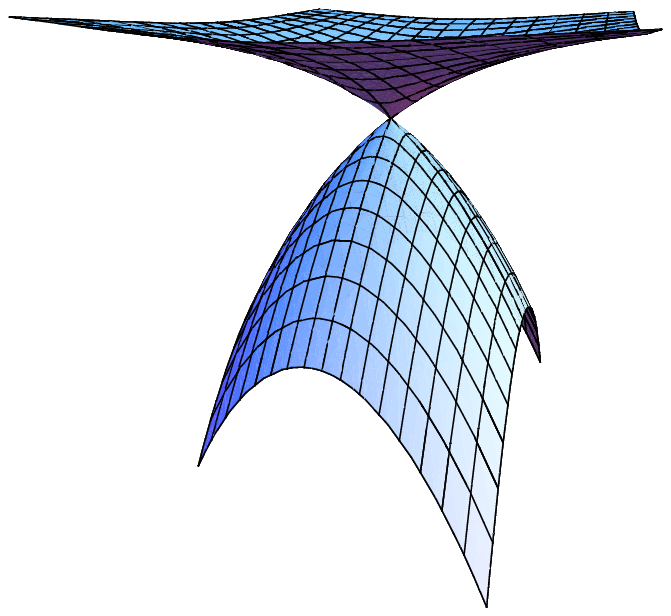}}
}
\caption{The positive and negative $\b$-exponential cones for $\b>0$ on $\R^2$, and the two shown together.}
\label{fig:cones}
\end{figure}

In $\R^n$ with the Euclidean metric $d(x,y)$, a $\db$-exponential cone~(\ref{expcone}) solves (\ref{biaslap}) away from $x_0$. In $\R$, the unique positive $\db$-exponential cone centered at $0$ and with values $C^+_0(0)=0$ and $C^+_0(1)=1$, namely $C^+_0(x)=\frac{1-e^{-\db|x|}}{1-e^{-\db}}$, coincides on $x\in[0,1]$ with $C^-_{1}(x)=\frac{e^\db-e^{\db |x-1|}}{e^\db-1}$, the negative $\db$-exponential cone that is centered at $1$ and has those values at 0 and 1. This function is easily seen to be the $\eps\to 0$ limit of the value function for the $\b$-biased $\eps$-game on $X=[0,1]$, $Y= \{ 0,1 \}$, with $F(0)=0$ and $F(1)=1$,  where the optimal strategies are obvious, and the game-play is just a biased random walk.

\begin{definition}\label{defCEC}
A function $u:X\setminus Y\longrightarrow\R$ satisfies {comparison with $\b$-exponential cones from above} (CECA) if for any $V\subset\subset X\setminus Y$ (i.e., $V$ is a bounded open subset of $X$ at a positive distance from $Y$), $x_0\in X$, $\iota \in\{+,-\}$ and $B\in\R$, whenever $u(x)\leq C^\iota_{x_0}(x) + B$ holds on $x\in\partial(V\setminus\{x_0\})$, we also have $u\leq C^\iota_{x_0}+B$ on the entire $V$. {Comparison with $\b$-exponential cones from below} (CECB) is defined analogously. If both conditions hold, then $u$ satisfies {comparison with $\b$-exponential cones} (CEC).
\end{definition}

The following equivalence is a generalization of the main observation of \cite{CEG}.

\begin{theorem}
\label{t.equiv}
Given an open connected $U\subset\R^n$, a continuous function $u: U\longrightarrow \mathbb{R}$ is a viscosity solution of (\ref{biaslap}) if{f} it satisfies comparison with $\db$-exponential cones.
\end{theorem}

Given this characterization of viscosity solutions, Theorem~\ref{t.main} is an immediate consequence of the following more general result.

\begin{theorem}
\label{t.general}
Let $X$ be a compact length space, with a specified closed subset $Y\subset X$ and $F$ a Lipschitz function on $Y$. Fix $\b \in \R$. Then there exists a unique extension of $F$ to $X$ that satisfies comparison with $\db$-exponential cones on $X\setminus Y$; it is the uniform limit of the value functions $u^\eps$ for the $\b$-biased $\eps$-tug-of-war games on $X$.
\end{theorem}

To prove convergence of the value functions $u^\eps(x)$, one would like to compare them for different values of $\eps$. It is unclear how to do this in general, but it works for $\eps$ and $2\eps$ as follows. One can view the $\eps$-game in rounds: until one player wins two more coin tosses than the other. The results of these rounds can then be considered as moves in the $2\eps$-game. The odds for player I to win a round is $\rho(\eps)^2$ to 1; this is because the coin toss trajectories ($T_k:=$ number of flips won by player I minus those won by player II in $k$ steps) that end with player I gaining the advantage of two flips are exactly the $T_k\mapsto -T_k$ reflections of the other trajectories, and each has probability $\rho(\eps)^2$ times that of the reflected one. So, if we want to mimic the $2\eps$-game exactly, we need $\rho(\eps)^2=\rho(2\eps)$. The simplest solution of this equation is $\rho_0(\eps)=e^{\beta\eps}$, which explains the special role of this odds function. More generally, if~(\ref{e.rhocond}) holds, we will be able to show monotonicity and hence pointwise convergence of the value functions along dyadic subsequences $\eps_n=\eps/2^n$, which will be the first step towards full convergence in the sup-norm.

Another manifestation of the special role of the odds function $\rho_0(\eps)=e^{\beta \eps}$ and the corresponding bias $\theta_0(\eps)$  is the following. If we perform a biased random walk $(X_n)_{n\geq 0}$ on $\R$ with steps $\pm\eps$, with bias $\theta_0(\eps)$ to the right, and $u(x):=a(1-e^{-\db x})$ is a one-dimensional $\b$-biased infinity harmonic function, then $u(X_n)$ is exactly a martingale.

We now define the {\bf value functions} $u_I^\eps$ for player I, and $u_{II}^\eps$ for player II. A strategy for a player is a way of choosing the player's next move as a function of all previous coin tosses and moves. Given two strategies, $\S_I$ for player I and $\S_{II}$ for player II, we define the payoff functions $F_+(\S_I,\S_{II})$ and $F_-(\S_I,\S_{II})$ as the expected payoff at the end of the game if the game terminates with probability one under these strategies. So, in this case they are equal. However, we let $F_-(\S_I,\S_{II})=-\infty$ and $F_+(\S_I,\S_{II})=\infty$ otherwise. We think of $u_I^\eps$ as the minimum that player I can guarantee being paid; on the other hand, player II can avoid having to pay more than the value $u_{II}^\eps$. More precisely,
\begin{eqnarray*}
  u_I^\eps & := & \sup_{\S_I}\inf_{\S_{II}}F_-(\S_I,\S_{II})\,,\\
  u_{II}^\eps & := & \inf_{\S_{II}}\sup_{\S_I}F_+(\S_I,\S_{II})\,.
\end{eqnarray*}
These values are functions of the starting location $x\in X$. Clearly, $u_I^\eps(x) \leq u_{II}^\eps(x)$. If they are equal, we say the game has a value,
and denote it by $u^\eps(x):=u_I^\eps(x)=u_{II}^\eps(x)$.
\medskip

The organization of the paper into sections is based on the following six main steps of the proof. For readers familiar with \cite{PSSW2}, we point out in brackets the main differences compared to the unbiased case.

\begin{enumerate}

\item $\theta(\eps)$-\textsc{infinity harmonicity}: The value functions
  $u_I^\eps(x)$ and $u_{II}^\eps(x)$ are
  $\theta(\eps)$-infinity harmonic (defined below).

\item \textsc{Game value existence}: The $\eps$-game
has a value, i.e., $u^\eps(x):=u_I^\eps(x)=u_{II}^\eps(x)$, for any $\theta(\eps)\in [-1,1]$. [To ensure finite stopping times, the backtracking argument now needs to be much more elaborate. In particular, the compactness of $X$ is used in a crucial way.]

\item \textsc{Convergence}: The pointwise limit $\vlim_{,\rho}:=\lim_{n\to\infty} u^{\eps/2^n}$ exists, provided condition~(\ref{e.rhocond}) holds. [The advantage given to one of the players becomes greater, hence harder to control, in one of the so-called ``favored'' games in Section~\ref{s.conv}.]

\item \textsc{Comparison with $\db$-exponential cones}:
  \begin{enumerate}
  \item The value function $u^\eps$ approximately
    satisfies comparison with $\db$-exponential cones. [The proof requires that sometimes the weaker player pulls towards a target; this is a reasonable strategy again because $X$ is compact.]
  \item The approximate CEC property provides enough regularity to conclude that, under~(\ref{e.rhocond}), $u^{\eps/2^n}$ converges not only pointwise, but also uniformly. [Improving the convergence from pointwise to uniform needs the regularity because we are allowing a wider class of biases, not just the nicest example $\rho_0(\eps)$.]
  \item Any uniform subsequential limit $\lim_{n\to\infty} u^{\eps_n} $ satisfies
CEC exactly.
  \end{enumerate}

\item \textsc{Uniqueness}: Any extension $u$ of $F$ that satisfies comparison with $\db$-exponential cones on $X\setminus Y$ is equal to the dyadic limit $\vlim_{,\rho_0}$. [Here it is important to consider the games with odds function exactly $\rho_0(\eps)$.] Hence all extensions $u$ and all limits $\vlim_{,\rho}$ are the same, independently of $\rho(\cdot)$ and $\eps>0$.

\item \textsc{Equivalence}: In the $\R^n$ case, a function satisfies CEC if and only if it is a viscosity solution to (\ref{biaslap}). [This generalizes \cite{CEG}.]

  \end{enumerate}
Finally, Section~\ref{s.open} collects some open problems. 
\medskip

In what follows, without loss of generality, we may assume that $\b>0$ and $\eps$ is small enough, hence the coin favors player I, and $\half < \frac{1+\theta(\eps)}{2} <1$.
\medskip

We now define $\theta(\eps)$-infinity harmonicity, the discrete analogue of having a constant drift in the direction of the gradient, and prove Item (1) above.

Starting from $x_0$, if player I wins the first coin toss and moves to some $x_1$, then conditioned on this event, he now has a strategy that guarantees a payoff arbitrarily close to $u_I^\eps(x_1)$, but does not have a strategy that guarantees more. The same is true if player II wins the first coin toss, and moves to some $x_1$ of his choice. The value $u_I^\eps(x_1)$ of the first choice is between the infimum and the supremum of $u_I^\eps$ over the ball $B_\eps(x_0)$, and can be arbitrarily close to these two extreme values. The same argument can be repeated for $u_{II}^\eps$, and we obtain:
\begin{lemma}
\label{l.infharm}
  For all $x$ in $X \backslash Y$, both $u_I^\eps$ and $u_{II}^\eps$ satisfy
\begin{equation}
  u(x)=\frac{1+\theta(\eps)}{2} \sup_{y\in B_\eps (x)} u(y) +
  \frac{1-\theta(\eps)}{2}\inf_{y\in B_\eps (x)} u(y) \,.
\label{e.infharm}
\end{equation}
\end{lemma}
\qed

We say that $u$ is {\bf $\theta(\eps)$-infinity harmonic} if it satisfies equation~(\ref{e.infharm}). Notice that the term $\eps$ in ``$\theta(\eps)$-infinity harmonic'' indicates both the bias and the size of the balls.

The supremum of the possible increases in the $\eps$-ball around $x$ is denoted by $\delta^+(x):=\sup_{y\in B_\eps (x)} u(y) -u(x)$. Similarly, the supremum of the possible decreases is $\delta^-(x):= u(x)- \inf_{y\in B_\eps (x)} u(y)$. Notice that, for a $\theta(\eps)$-infinity harmonic function $u$, equation (\ref{e.infharm}) implies that these values are related by the equation
\begin{equation}
\delta^+(x) =
  \frac{1-\theta(\eps)}{1+\theta(\eps)}\, \delta^-(x) = \frac{\delta^-(x)}{\rho(\eps)} \,,
\label{e.deltarel}
\end{equation}
and we see that $\delta^+(x) < \delta^-(x)$.

We will sometimes use the modified metric $d^\eps(x,y)$, which is $\eps$ times the minimum number of steps of size less than $\eps$ to get from $x$ to $y$; for $x\not= y$ this is $\eps+\eps\lfloor d(x,y)/\eps\rfloor$. Let us remark here that we use open $\eps$-balls for the steps because, in Section~\ref{s.conv} , we will need that any $2\eps$-step can be broken into two $\eps$-steps without touching $Y$ in between.

Whenever we write $f=O(g)$ for two functions, we mean that $f\leq Cg$, where the constant $C\in (0,\infty)$ might depend on $\beta$ and $X,Y,F$.

%
%
\section{Game value existence}\label{s.value}

We will always assume that $X$ is a compact length space, $Y\subset X$ is closed, and $F:Y\longrightarrow \R$ is Lipschitz.

\begin{lemma}
$u_I^\eps$ is the smallest $\theta(\eps)$-infinity harmonic function bounded from
below on $X$ that extends $F$. More generally, if $v$ is a $\theta(\eps)$-harmonic function which is bounded from below on $X$ and $v \geq F$ on $Y$, then $v \geq u_I^\eps$ on $X$.

Similarly, $u_{II}^\eps$ is the largest $\theta(\eps)$-infinity harmonic
function which is bounded from above on $X$ and extends $F$.
\label{small}
\end{lemma}

\begin{proof}
  From Lemma \ref{l.infharm} we know that $u_I^\eps$ is $\theta(\eps)$-infinity harmonic. It is bounded
  below because player I can fix a $y \in Y$ and adopt a ``pull toward
  $y$'' strategy which at each move attempts to minimize the distance of the game position to $y$. In an infinite series of coin tosses there will
  be a time when player I has won more tosses by the amount $d^\eps(x_0,y)/\eps$,
  therefore he can ensure that the game terminates. Hence
  $u_I^\eps(x_0) \geq \inf_YF > -\infty$. Of course $\inf_YF > -\infty$ because $F$ is a Lipschitz function and $Y$ is compact.

  Suppose $v \geq F$ on $Y$. Suppose also that $v$ is $\theta(\eps)$-harmonic
 and bounded below on $X$. In order to show that $v \geq
  u_I^\eps =\sup_{\S_I} \inf_{\S_{II}} F_-(\S_I,\S_{II})$, we must show
  that for every $\S_I$ there exists an $\S_{II}$ such that $v \geq
  F_-(\S_I,\S_{II})$. We may assume that $\S_I$
  is a strategy that ensures that the game terminates --- otherwise $F_-(\S_I,\S_{II})=-\infty$ for some $\S_{II}$ and the assertion is obvious. We select $\S_{II}$
  such that if player II wins the $k^\text{th}$ coin toss he moves to $x_k$ where
  $v(x_k) < \inf_{y\in B_\eps(x_{k-1})}v(y) + \frac{\delta}{2^k}$. In other words, he plays to ``almost'' minimize $v$. 
  Now notice that $M_k=v(x_k)+ \delta 2^{-k}$ is a supermartingale with respect to the filtration $\cal F_k$ generated by the first $k$ coin tosses and moves:
  \begin{eqnarray*}
    \mathbb{E}[M_k | {\cal{F}}_{k-1}] & \leq & \frac{1+\theta(\eps)}{2}
    \sup_{y \in B_\eps (x_{k-1})}v(y) + \frac{1-\theta(\eps)}{2}
    \left(\inf_{y\in B_\eps (x_{k-1})}v(y) + \frac{\delta}{2^k}\right)
    +\frac{\delta}{2^k}\\
    &=& v(x_{k-1}) +\left(\frac{1-\theta(\eps)}{2} +1\right) \frac{\delta}{2^k}\\
    &<& v(x_{k-1}) +2 \frac{\delta}{2^k} = v(x_{k-1})+ \frac{\delta}{2^{k-1}} = M_{k-1}\,,
  \end{eqnarray*}
  where the equality follows because $v$ is $\theta(\eps)$-harmonic.

  Let $\tau:= \inf \{k : x_k \in Y\}$; we have assumed that $\S_I$ an $\S_{II}$ are such that  $\tau<\infty$ a.s. Since $M_k$ is a supermartingale bounded from below,
an optional sampling theorem (proved by monotone convergence and Fatou's lemma) gives $ \mathbb{E} [M_\tau] \leq \mathbb{E}[M_0]$ and we have
  $$
  F_-(\S_I,\S_{II}) = \mathbb{E}[F(x_\tau)] \leq \mathbb{E}[v(x_\tau)] \leq \mathbb{E} [M_\tau] \leq
  \mathbb{E}[M_0] = v(x_0) + \delta,
  $$
  where $\delta$ was arbitrary, hence $u_I^\eps(x_0) \leq v(x_0)$.
\end{proof}

\begin{theorem}
  For $X,Y,F$ as before, for any $\eps>0$ and $\theta(\eps)\in [-1,1]$, we have $u_I^\eps=u_{II}^\eps$, i.e., the $\eps$-game has a value.
\label{gameval}
\end{theorem}

\begin{proof}
  We know that $u_I^\eps \leq u_{II}^\eps$ always holds, so we need to
  show $u_I^\eps \geq u_{II}^\eps$. As in Lemma \ref{small}, we
  know that $u_{II}^\eps \leq \sup_YF < \infty$. We let
  $u=u_{II}^\eps$. Recall that $\delta^+(x)= \sup_{y \in
    B_\eps(x)}u(y)-u(x)$; we will use the shorthand notation $\delta_k:=\delta^+(x_k)$. We first
  assume that the infima and suprema in the definition of $\theta(\eps)$-harmonicity
 are achieved for all $x\in X$, and that $\delta_0 >0$. Always achieving the extrema is of course unrealistic, especially since the balls $B_\eps(x)$ are open. However, this way the formulas will be simpler and the main ideas clearer --- later we will remove these assumptions.

For $m\in\Z^+\cup\{\infty\}$, a finite or infinite sequence $(x_i : 0 \leq i < m)$ of possible $\eps$-game positions in $X\setminus Y$ is called a maximal chain if the following holds: $u(x_{k+1})=u(x_k)+\delta_k$ for $k<m$, and if $m<\infty$, then the chain cannot be continued, i.e., all points $x_m\in X$ with $u(x_m)=u(x_{m-1})+\delta_{m-1}$ are in fact in $Y$. Notice that $u(x_{k+2}) > u(x_{k+1})$ and, since the
  maximal value in $B_\eps(x_k)$ was attained at $x_{k+1}$,
  the point $x_{k+2}$ cannot be in $ B_\eps (x_k)$. Similarly, $x_j$ is not in
  $B_\eps (x_k)$ for any $j \geq k+2$. So, for every $k\ne j$ we have
  $d(x_{2k},x_{2j}) \geq \eps$, i.e., the open balls $B_{\eps/2}(x_{2k})$ are all disjoint.
Now, since $X$ is compact, there exists a finite covering of it by open $\eps/4$-balls. Each ball $B_{\eps/2}(x_{2k})$ contains one of these $\eps/4$-balls entirely, hence there are more $\eps/4$-balls in the covering than points $x_{2k}$. Therefore, there is a uniform finite upper bound $M=M(\eps)$ on the cardinality of any maximal chain $\{ x_0,x_1,x_2,...,x_{m-1}\}$, which is independent of the starting point $x_0$.

  We now choose
  $\delta$ small enough so that $\rho^M \delta < \delta_0$, where $\rho=\rho(\eps)$.  We
  let $A_0:=Y$ and for $n>0$ we let $A_n:=\big\{x \in X : \exists z\in
  A_{n-1}$ with $d(z,x)<\eps$ and $u(z)=\sup_{y \in B_\eps(x)}u(y) >u(x)+
  \rho^{n-1} \delta \big\}$. Finally, $A:= \bigcup_{n=0}^M A_n$.

  Notice that $x_0$ is in $A$. For, if $\{x_0,x_1,x_2,...,x_{m-1} \}$ is
  the maximal chain described above, then its length is at most $M$.
  The $\theta(\eps)$-harmonicity condition for $u$ and (\ref{e.deltarel}) imply inductively that $\delta_{m-1} \geq \rho^{-m+1} \delta_0 \geq
  \rho^{-M+1} \delta_0 >\delta$, so $x_{m-1}$ is in
  $A_1$ and hence $x_0$ is in $A$. Similarly, every $x\not\in A$ satisfies $\delta^+(x)\leq \rho^{M-1} \delta$.

  For each $n$, we let $j_n := \max \{ j \leq n : x_j \in A\}$
  and $v_n:=x_{j_n}$; this is the last position in $A$ up to time $n$.

  We need to show that there exists an $\S_I$ such that for every
  $\S_{II}$, we have  $F_-(\S_I,\S_{II})\geq u$. To that end, we let
  $\S_{II}$ be an arbitrary strategy and we let $\S_I$ be the following
  strategy for player I. When $x_n$ is not in $A$ and player I wins
  the coin toss he backtracks toward $A$ by moving to a previous game
  position that is closer to $v_n$. That is, he moves to any position
  in $\{ x_{j_n}, x_{j_n+1}, ... , x_{n-1} \}$ which is within
  $\eps$ of $x_n$ and is closer to $v_n$ than $x_n$. When $x_n$ is
  in $A$, he maximizes $u$ within $A$ as much as possible.

  We let $d_n$ be the minimum number of steps of length less than $\eps$
  required to get from $x_n$ to $v_n$ by only moving to
  previous game positions.  We define
  \be\label{subMG}
  m_n:=u(v_n) +(n-d_n)\,\rho^M\delta\,.
  \ee
  We claim that $m_n$ is a submartingale.

  Case I. We consider the case where $x_n \in A$ and so $m_n=u(x_n) + n\rho^{M}\delta$.

  1. If player I wins the coin toss, then his strategy determines that
  $m_{n+1} = u(x_{n+1}) + (n+1)\rho^{M}\delta=u(x_n)+\delta_n +
  (n+1)\rho^{M}\delta > m_n+\delta_n$.

  2. a. If II wins the coin toss and moves to another position
  $x_{n+1} \in A$, then $m_{n+1} = u(x_{n+1}) +
  (n+1)\rho^{M}\delta \geq u(x_n)- \rho \delta_n +
  (n+1)\rho^{M}\delta > m_n- \rho \delta_n$.

  b. If he moves to a position outside of $A$, then $m_{n+1} = u(x_n) +
  ((n+1)-1) \rho^{M}\delta= m_n$.

  Thus, in this case $m_n$ is a submartingale, since
  $$
  \mathbb{E}[m_{n+1}|{\cal{F}}_n] \geq \frac{1+\theta(\eps)}{2} (m_n +
  \delta_n) +\frac{1-\theta(\eps)}{2} \left(m_n - \epsfrac \delta_n\right)= m_n\,.
  $$

  Case II. We next consider the case where $x_n$ is not in $A$ and $m_n=u(v_n) + (n-d_n)\rho^{M} \delta$, with $d_n >0$.

  1. If player I wins the coin toss, then his strategy determines that he
  moves at least one $\eps$ step closer to $v_n$. Hence $m_{n+1}
  \geq u(v_n) + (n+1-(d_n-1)) \rho^{M} \delta = m_n + 2
  \rho^{M} \delta$.

  2. a. If II wins the coin toss and moves to another position not in
  $A$, we have $m_{n+1} \geq u(v_n) + ((n+1) -(d_n+1))
  \rho^{M}\delta = m_n$.

  b. If II wins and moves back into $A$, then we have
\begin{eqnarray*} m_{n+1}
  &=&u(x_{n+1}) + (n+1) \rho^{M} \delta \\ 
  &\geq& u(x_n) -\rho \delta_n + (n+1) \rho^{M} \delta \\
  &\geq& u(x_n) + n \rho^{M} \delta 
  \geq u(v_n) +(n-d_n) \rho^{M} \delta
  = m_n\,;
\end{eqnarray*}
from the second to the third line we used that $\delta_n \leq \rho^{M-1}\delta$, since $x_n\not\in A$, while in the third line we used that there is a path between $x_n$ and $v_n$ with at most $d_n$ steps, each outside $A$, hence having a value difference at most $\rho^{M} \delta$.

  Altogether, in Case II we also have
  $$
  \mathbb{E}[m_{n+1}|{\cal{F}}_n] \geq m_n\,.
  $$

  Hence $m_n$ is indeed a submartingale. We let $\tau$ be the first
  time a terminal state is reached. We have that for all $n$,
  $\mathbb{E}[m_{n\wedge \tau}] \geq \mathbb{E}[m_0] = u(x_0)$.

  We now show that $\mathbb{E}\tau$ is finite, so that $\mathbb{E}F(x_\tau)\geq u(x_0)$.

  Whenever the game position is in $A$, there is an independent chance
  that player I will win the necessary tosses in a row (some number
  at most M) and end the game. So the number of steps in $A$ is
  stochastically bounded by $M\mu$ where $\mu$ is a geometric random
  variable with success probability $(\frac{1+\theta(\eps)}{2})^M$. Consequently, the
  expected number of steps inside $A$ is at most
  $\mathbb{E}[M\mu]=M(\frac{2}{1+\theta(\eps)})^M$.

  Each outside run is stochastically bounded by an independent random
  variable $\xi_k$, the length of a random walk on the nonnegative
  integers started from 1 until hitting 0, where the probability of moving to the left is
  $\frac{1+\theta(\eps)}{2}$ and to the right is
  $\frac{1-\theta(\eps)}{2}$. We know that
  $\mathbb{E}\xi_k=\frac{1}{\theta(\eps)}$.

  The number of runs outside $A$ is bounded by the number of steps in
  $A$. This gives $$\mathbb{E}\tau \leq \mathbb{E}(\sum_{k\leq M\mu}
  \xi_k + M\mu).$$ By Wald's Lemma we have that $\mathbb{E}
  \sum_{k\leq M\mu} \xi_k \leq M (\frac{2}{1+\theta(\eps)})^M
  \frac{1}{\theta(\eps)}$, and so
$$\mathbb{E}\tau \leq M
  \left(\frac{2}{1+\theta(\eps)}\right)^M \left(\frac{1}{\theta(\eps)}+1\right)=:T\,.$$
So $\mathbb{E}m_\tau \geq \mathbb{E}m_0=u(x_0)$.

  Recall that we have seen that $u(x_n) \geq m_n -n \rho^{M}
  \delta$ for all $x_n$, by considering the path between $x_n$ and $v_n$ outside $A$.

  Summarizing, we have that
\be\label{Tfinal}
\mathbb{E}F(x_\tau ) \geq \mathbb{E}m_\tau -
  \rho^{M} \delta\, \mathbb{E}\tau \geq m_0- \rho^{M} \delta
  T = u(x_0)- \rho^{M} \delta T\, .
\ee
Since $\delta$ can be made
  arbitrarily small, we have $\mathbb{E}F(x_\tau) \geq
  u(x_0)=u_{II}^\eps(x_0)$. Thus $u_I^\eps \geq
  u_{II}^\eps$, and, in the case where the suprema and infima are
  achieved, we are done.
\medskip

We now remove the assumption that the extrema are achieved, but still assume that $\delta_0=\delta^+(x_0)>0$. Recall that there is a uniform bound $M=M(\eps)$ on the length of any sequence $\{x_0,x_1,\dots,x_{m-1}\}$ where $d(x_j,x_k) \geq \eps$ for $j \geq k+2$. Let $\delta,\eta>0$ be small reals such that
\be\label{e.deltaeta}
\delta > (\rho+2)\eta\qquad\text{and}\qquad \delta_0 > \rho^M(\delta+M\eta)\,.
\ee
Now, given $x_k$, choose $x_{k+1}$ to be a point in $B_\eps(x_k)$ with
\be\label{k+1}
u(x_{k+1}) > u(x_k) + \delta_k - \rho^{M-k}\eta\,.
\ee
We claim that $u(x_{k+1})\geq u(x_k)$ and $u(x_{k+2}) > u(x_k) + \delta_k$, therefore $d(x_j,x_k) \geq \eps$ for $j \geq k+2$, and the sequence must reach $Y$ in at most $M$ steps.

By (\ref{e.deltarel}) and (\ref{k+1}), we have $\delta_{k+1} \geq \rho^{-1}(u(x_{k+1})-u(x_k)) > \rho^{-1} (\delta_k - \rho^{M-k}\eta)$. This and the second inequality of (\ref{e.deltaeta}) give by induction that
\be\label{deltak}
\delta_k > \rho^{M-k}(\delta + (M-k)\eta)\,.
\ee
Combining this and (\ref{k+1}) gives $u(x_{k+1})\geq u(x_k)$ immediately. To show  $u(x_{k+2}) > u(x_k) + \delta_k$, we use (\ref{k+1}) twice to write $u(x_{k+2})$ in terms of $u(x_k)$, then using (\ref{deltak}) for $\delta_{k+1}$ we arrive at the following inequality to verify:
$\rho^{M-k-1}\big(\delta + (M-k-2-\rho)\eta\big) > 0$, which always holds by the first inequality of (\ref{e.deltaeta}), and our claim is proved.

This argument shows that if we now define the set $A:=\bigcup_{n=0}^M A_n$ by $A_0:=Y$ and
\begin{align*}
A_n:=\big\{x \in X &: \delta^+(x) > \rho^{n-1}\big(\delta + (n-1)\eta\big),\text{ and}\\
&\ \exists z\in A_{n-1}\text{ with }d(z,x)<\eps\text{ and }u(z)> u(x)+\delta^+(x)-\rho^{n-1}\eta \big\},
\end{align*}
then $x_0 \in A$, and similarly, $x\not\in A$ implies $\delta^+(x) \leq \rho^M(\delta + M\eta)$.

The strategy for player I now will be the following: if $x_n\in A$, say $x_n\in A_i$, then try to step to a point $z\in A_{i-1}$ with $u(z)>u(x_n)+\delta_n-\rho^{i-1}\eta$; if $x_n\not\in A$, then try to backtrack towards the last visit $v_n$ to $A$, as before. We also slightly modify the definition (\ref{subMG}), and claim that now
$$
m_n:=u(v_n)+(n-d_n)\,\rho^M(\delta+M\eta)
$$
is a submartingale. The proof remains the same, with two small modifications. Firstly, in Case I.1, we have only $u(x_{n+1}) > u(x_n) + \delta_n - \rho^M\eta$, but then $(n+1)\rho^M(\delta+M\eta)-\rho^M\eta > n\rho^M(\delta+M\eta)$ still gives us what we need. Secondly, in Case II.2.b, the path between $x_n$ and $v_n$ now makes possible a value difference at most $d_n \rho^M(\delta+M\eta)$, but we changed the definition of $m_n$ exactly in order to accommodate this difference.

The proof is completed by noting that given $\eps>0$, the term $\rho^M(\delta+M\eta)T$ in the present version of (\ref{Tfinal}) can be made arbitrarily small by choosing $\delta$ and $\eta$ small.
\medskip

Finally, if $\delta_0=0$, we fix $y \in Y$ and let player I pull toward $y$
  until a position $x_0^*$ with $\delta^+(x_0^*) >0$ is reached or the
  game terminates. Note that the first such $x_0^*$ encountered has $u(x^*_0)=u(x_0)$. If $x_0^* \in Y$, then $F(x_0^*)=u(x_0)$, and player I achieved the payoff $u(x_0)$. If $x_0^* \in X
  \backslash Y$, then player I plays the above strategy with $x_0^*$ replacing $x_0$, and can achieve the payoff $u(x^*_0)=u(x_0)$ in expectation. This concludes the proof.
\end{proof}

\begin{corollary}
  $u^\eps=u_I^\eps=u_{II}^\eps$ is the unique bounded
$\theta(\eps)$-infinity harmonic function agreeing with $F$ on $Y$.
\end{corollary}

\begin{proof}
  By Lemma \ref{small} we have that $u_I^\eps$ is the smallest such function and $u_{II}^\eps$
  is the largest. By Theorem \ref{gameval} they are equal. Thus the function is
  unique.
\end{proof}

%
%
\section{$u^{\eps/2^n}$ Convergence}\label{s.conv}

We seek to prove that $u^{\eps/2^n}$ converges as $n \in \mathbb{N}$ tends to $\infty$. However, it is not the case that $u^{\eps/2^n}(x)$ is monotone for each $x\in X$, so, as in \cite{PSSW2}, we define a new game that does possess this property, the {\em biased II-favored $\eps$-tug-of-war}. In this game, player II has some advantages: if he wins the coin toss, he is allowed to choose the next game position in a ball larger than the usual one; otherwise he still has the option to slightly modify the move of player I, provided that his modification terminates the game.

More precisely, the {\em biased II-favored $\eps$-tug-of-war} is defined as follows. At the $k^\mathrm{th}$ step, player I chooses $z_k$ in $B_\eps(x_{k-1})$. If player I wins the coin toss, then player II can choose $x_k \in (B_{2\eps}(z_k) \cap Y) \cup \{ z_k \}$. In other words, he moves to
$z_k$ or terminates the game. If player II wins the coin toss he moves to a position in $B_{2\eps}(z_k)$ of his choice.

The value to player I of the biased II-favored $\eps$ tug-of-war is denoted by $v^\eps$. The value to player II of the correspondingly defined
{\em biased I-favored $\eps$ tug-of-war game} is denoted by $w^\eps$.
In both cases, the bias is $\frac{1+\theta(\eps)}{2}$, favoring player I.

As the names suggest, player I cannot do better in the II-favored game than in the ordinary game. In other words, $v^\eps \leq u_I^\eps$. This is because given any strategy for player II in the ordinary $\eps$-game, player II can follow it also in the II-favored game, since for any $z\in B_\eps(x)$ we have $B_\eps(x) \subseteq B_{2\eps}(z)$.
Similarly, player II cannot do better in the I-favored game than in the ordinary game (he cannot ensure a value that is smaller than the value in the ordinary game), so $u_{II}^\eps\leq w^\eps$. Hence we have $v^\eps\leq u_{I}^\eps=u_{II}^\eps \leq w^\eps$. Moreover, the following lemma holds:

\begin{lemma}
\label{monot}
If $\rho(\cdot)$ is log-concave, then, for any $\eps>0$,
$$ v^{2\eps} \leq v^{\eps} \leq u^\eps \leq w^\eps.$$
If $\rho(\cdot)$ is log-convex, then
$$ v^{\eps} \leq u^\eps \leq w^\eps \leq w^{2\eps}.$$
\end{lemma}

\begin{proof}
The inequalities $v^{\eps} \leq u^\eps \leq w^\eps$ were already shown above.
We now prove $v^{2\eps} \leq v^\eps$ for the case when $\rho(\cdot)$ is log-concave. The argument for $w^{\eps} \leq w^{2\eps}$ is analogous.

Given any strategy $\S_I^{2\eps}$ for player I in the II-favored $2\eps$-game, he can mimic it in the II-favored $\eps$-game as follows. Whenever he would choose a target point $z\in B_{2\eps}(x)$ in $\S_I^{2\eps}$, in the $\eps$-game he aims towards $z$ until what we will call ``a round'': until one player wins two more coin tosses than the other. This person will be player I with
odds $\rho(\eps)^2$, as we explained in the Introduction (or as follows from the computation in Lemma~\ref{wbound} below).
Note that log-concavity implies that $\rho(\eps)^2 \geq \rho(2\eps)$, hence the odds for player I to win a round is larger than his odds for a single coin toss in the $2\eps$-game.

Now, if player I wins the round, the game position will be $z$, or in $Y\cap B_{4\eps}(z)$ if player II ended the game during this round; otherwise, it will be within $B_{4\eps}(z)$, as is easily checked. Therefore, player II has at most as many choices if he wins the round as he had upon winning a coin toss in the $2\eps$-game. Altogether, player I can do at least as well with this strategy as in the $2\eps$-game, and we are done.
\end{proof}

Thus, if $\rho(\cdot)$ is log-concave, then $(v^{\eps/2^n})_n$ is monotone increasing, hence its pointwise limit exists, while if $\rho(\cdot)$ is log-convex, then $\lim_{n \rightarrow \infty} w^{\eps/2^n}$ exists. However, a priori these limits could be $+\infty$ or $-\infty$ sometimes, may not be the same even if both exist, and could depend on $\rho(\cdot)$ and $\eps$.

\begin{lemma}
  Let $\eps >0$. Then for each $x \in X$ and $y \in Y$,
  $$
  v^\eps (x) \geq F(y)-(2\eps + d^\eps (x,y))
  \mathrm{Lip}_Y F.
  $$
  Such an expected payoff is guaranteed for player I if he adopts a
  ``pull towards $y$'' strategy which at each move attempts to reduce
  $d^\eps(x_k,y)$.
\label{vbound}
\end{lemma}

\begin{proof}
  Given $x=x_0 \in X$ and $y \in Y$, the ``pull towards $y$'' strategy for player I in the II-favored game ensures that the game ends almost surely in finite time $\tau<\infty$, and that $d^\eps(x_n,y)$ is a supermartingale, except at the
  last step where player II may have moved the game position up to
  $\eps$ farther from $y$ even if player I won the coin toss.
  This implies
  $$
  \mathbb{E}[d^\eps (x_\tau ,y)] < d^\eps (x,y) +2\eps.
  $$
  Since $F$ is Lipschitz on $Y$, we get
  $$
  \mathbb{E}[F(x_\tau )] \geq F(y) - (d^\eps(x,y) +2\eps)
  \mathrm{Lip}_Y F\,.
  $$
  Since $v^\eps$ is at least as great as the expected payoff under
  this pull toward $y$ strategy, the lemma is proved.
\end{proof}

We define $d^\eps (X)= \sup_{x_1,x_2 \in X}d^\eps(x_1,x_2)$. Since $X$ has finite diameter ($\mathrm{diam}\, X$) in the path metric, we have $d^\eps (X) <\infty$; in fact, $\lim_{\eps\to 0} d^\eps(X)=\mathrm{diam}\, X$.

\begin{lemma}
  If $x \in X$ , $y \in Y$ and $d^\eps(x,y)=\eps$, then
$$u^\eps(x) \leq F(y) +
  \rho(\eps)^{d^\eps(X)/\eps} \, \eps \, \mathrm{Lip}_Y F \leq
  F(y) + O(1)\, e^{\b\,\mathrm{diam}\, X} \, \eps \, \mathrm{Lip}_Y F,
$$
with an absolute constant $O(1)$.
\label{ubound}
\end{lemma}

\begin{proof}
  Let $u=u_I^\eps$.  Let $\gamma= u(x) - F(y) >0$. (If $\gamma \leq 0$, we are done).
  
As in the proof of Theorem~\ref{gameval}, we first assume that the infima and suprema in the definition of $\theta(\eps)$-harmonicity are always achieved. Then there is a path $\{x=x_0,x_1,x_2,...\}$ of maximal
ascent, i.e., for each $k$, we have $u(x_k)=\sup_{z\in B_\eps(x_{k-1})} u(z)$.
Since $X$ is bounded, this sequence is finite and ends at some $y' \in Y$.

Let $d$ be the minimal number of steps of size less than $\eps$ required to get from $x$ to $y'$, i.e., $d=d^\eps(x,y')/\eps$. So $y'=x_n$ for some $n \geq d$.

  We let $\delta_k= u(x_{k+1})-u(x_k)$ be the maximal possible
  increase of $u$ in the $\eps$-ball around $x_k$. Recall that at
  each step, $\delta_k \geq \delta_{k-1}/\rho(\eps)$ and $\delta_0 \geq
 \gamma/\rho(\eps)$. So $\delta_k \geq
 \rho(\eps)^{-(k+1)} \gamma$, hence
  $$
  F(y') \geq u(x_0) + \sum_{k=0}^{n-1}
  \rho(\eps)^{-(k+1)} \gamma \geq u(x_0) +
  \sum_{k=1}^d \rho(\eps)^{-k} \gamma.
  $$
  On the other hand, $d(y,y') \leq d^\eps(y,y') \leq (d+1)\eps$, and since $F$
  is Lipschitz, we have
  $$
  F(y') \leq F(y) + (d+1) \eps \mathrm{Lip}_Y F.
  $$
  Combining these two inequalities and using $u(x_0)-F(y) = \gamma$, we get
  $$
  (d+1) \eps \mathrm{Lip}_Y F \geq \sum_{k=0}^d
  \rho(\eps)^{-k} \gamma \geq (d+1)
  \rho(\eps)^{-d} \gamma \geq (d+1)
  \rho(\eps)^{-d^\eps (X)/\eps} \gamma,
  $$
  and so
  $$
  \rho(\eps)^{d^\eps (X)/\eps} \, \eps \, \mathrm{Lip}_Y F \geq \gamma.
  $$
  Thus, in the case where the infima and suprema are achieved, we are
  done.

  If the infima and suprema are not achieved, we fix $\eta,\delta > 0$ small enough for~(\ref{e.deltaeta}) to hold, with our $M=M(\eps)$ as before. (Note here that we had $\gamma>0$, hence also $\delta_0>0$.) We now define our sequence
  $\{x_0,x_1,x_2,...\}$ such that each $x_k$ satisfies~(\ref{k+1}). Similarly to what we did after~(\ref{k+1}), we have 
 $\delta_{k+1}\geq \rho^{-1}(\delta_k-\rho^{M-k}\eta)$, which by induction gives $\delta_{k+1} \geq \rho^{-(k+1)}\delta_0-(k+1) \rho^{M-k-1}\eta$. Since $\delta_0 \geq \gamma/\rho$, we get
\begin{align*}
  F(y') \geq u(x_0) + \sum_{k=0}^{d-1} \left(\delta_k-\rho^{M-k}\eta\right) 
 & \geq u(x_0) + \sum_{k=0}^{d-1} \left( \rho^{-(k+1)} \gamma - (k+1)\rho^{M-k}\eta \right) \\
 & \geq u(x_0) + \sum_{k=1}^d \rho^{-k} \gamma - d^2\rho^M \eta\,.
\end{align*}
Since $d\leq M$, from this we get in the same way as above that 
  $$
  \rho(\eps)^{d^\eps (X)/\eps} \, \eps \, \mathrm{Lip}_Y F \geq \gamma - M \rho^M \eta\,.
  $$
Since $\eta$ was arbitrary, the proof is complete.
\end{proof}

Some version of the above lemma is probably valid also for $w^\eps$ in place of $u^\eps$, but we prove instead a weaker but more general bound for $w^\eps$, not using the boundedness of $X$.

\begin{lemma}
  If $x \in X$ and $y \in Y$ with $d^\eps(x,y)=\eta \in [\eps,1)$, then, for
  small enough $\eps$, we have $w^\eps(x) \leq F(y) +
  \sqrt{\eta} \, (\mathrm{Lip}_Y F + \exp(4\b\sqrt{\eta})\sup_{y' \in
    Y}F(y'))$.
\label{wbound}
\end{lemma}

\begin{proof}
  In the I-favored game, let player II pull toward $y$. The probability that
  player I gets $\sqrt{\eta}$ distance away from $x$ before the game ends is at most $O(\sqrt{\eta})$, which can be seen as follows.  Let  $k=\lfloor \sqrt{\eta}/{\eps} \rfloor-1$ and $\ell=\eta/\eps \in \Z_+$. Let $R_0=0$, and $R_{-j}=-\sum_{i=1}^{j} \rho^{-i}$ for $0<j\leq k$, and $R_j=\sum_{i=0}^{j-1} \rho^i$ for $0<j\leq \ell$. Consider the one-dimensional game that starts at $R_0$, player I moves to the left with
  probability $\cefrac$, player II moves to the right with
  probability $\fracce$, and the end is $\tau=\inf \{n: x_n=R_{-k} \text{ or }
  x_n=R_{\ell}\}$. Let $q$ be the probability that
  $x_\tau=R_{-k}$. This upper bounds the probability that player I will
  get $\sqrt{\eta} \geq (k+1)\eps$ distance away from $x$ before the original (I-favored) game ends. 
  The position in the one-dimensional game is a
  martingale, so $0=\mathbb{E}x_0=\mathbb{E} x_\tau=q R_{-k} + (1-q) R_{\ell}$. Hence
  $$
  q=
  \frac{R_\ell}{R_\ell-R_{-k}} < \frac{\ell\rho^{\ell-1}}{(k+\ell)\rho^{-k}} < \frac{\eta/\eps}{\sqrt{\eta}/\eps-1}\,\rho(\eps)^{\sqrt{\eta}/\eps+\eta/\eps} < \sqrt{\eta} \exp(4\b\sqrt{\eta})\,,
  $$
where the last inequality uses that $\rho(\eps)<1+2\b\eps$ for $\eps$ small. 

Now, since all points in $Y$ which are less than $k$ steps away from $x$
  have value at most $F(y)+ k\eps \, \mathrm{Lip}_Y F$, we have
  $$
  w_\eps (x) \leq F(y) + \sqrt{\eta}\, \mathrm{Lip}_Y F +
  \sqrt{\eta} \exp(4\b\sqrt{\eta}) \sup_{y' \in Y} F(y')\,,
  $$
and we are done.
\end{proof}

\begin{proposition}
  $\|v^\eps - u^\eps\|_\infty = O(\eps)$.
\label{vconv}
\end{proposition}

\begin{proof} We will closely follow the proof of \cite[Lemma 3.3]{PSSW2}.

 We play a virtual ordinary game, and player I's strategy in the
  favored game will mimic his strategy in the virtual game so that the
  game positions stay close (within $\eps$ of each other). More precisely, we
  first fix a strategy $\S_I$ for player I in the ordinary $\eps$-game
  that guarantees an expected payoff of at least $u_I^\eps -\eps$ against any strategy of player II. (In other words, $F_-(\S_I,\S_{II}) \geq u_I^\eps - \eps$ for any $\S_{II}$.) Then, we will select a strategy $\S_I^{\cal{F}}$ for player I in the favored game which does almost as well as $\S_I$, in the sense that for any strategy $\S_{II}^{\cal{F}}$ for player II, there will exist an $\S_{II}$
such that $F_-(\S_I^{\cal{F}},\S_{II}^{\cal{F}}) \geq F_-(\S_I,\S_{II})- O(\eps) \geq u_I^\eps - O(\eps)$, where the last inequality is by our assumption on $\S_I$. Thus we will have $v^\eps \geq u_I^\eps - O(\eps)$. Combining this with $u_I^\eps \geq v^\eps$ from Lemma~\ref{monot}, we will have the desired result.

  We let the superscript $\cal{F}$ denote positions in the favored
  game. We begin both games at the same position. In other words, we
  let $x_0^{\cal{F}}=x_0$. Our strategy $\S_I^{\cal{F}}$ will satisfy
  by induction that $d(x_n^{\cal{F}},x_n) < \eps$ for each $n$, as long as both games are running.

  Assuming this for some $n$, at step $n+1$ player I chooses
  $z_{n+1}^{\cal{F}}=x_n$. In other words, he aims for the current state in
  the virtual game.

  If player I wins the coin toss, then
  $x_{n+1}^{\cal{F}}=z_{n+1}^{\cal{F}}=x_n$ (as long as the game has not
  terminated), and in the virtual game $x_{n+1}$ is chosen according to
  $\S_I$. Hence $d (x_{n+1} , x_n) < \eps$, and thus $d
  (x_{n+1}^{\cal{F}} , x_{n+1}) < \eps$.

  If player II wins the coin toss, then he chooses some $x_{n+1}^{\cal{F}}$ according to his $\S_{II}^\cal{F}$, satisfying $d (x_{n+1}^{\cal{F}},
  x_{n}) < 2\eps$. Then, in the strategy $\S_{II}$ we are constructing, player II
  picks $x_{n+1}$ such that $d (x_{n+1}^{\cal{F}} , x_{n+1}) <
  \eps$ and $d (x_n , x_{n+1}) < \eps$.

  So, the game positions stay ``close'' in either case.

If the virtual ordinary game terminates before the favored game does, in $\S_I^\cal{F}$ we let player I pull towards the terminal position of the virtual game. If the favored game terminates before the virtual one, then in $\S_{II}$ we let player II pull towards the final position of the favored game. We let $\tau$ be the time the virtual game ends and let $\tau^{\cal{F}}$ be the time when the favored game ends.

  Case I. If $\tau=\tau^{\cal{F}}$, then $x_\tau$ and $x_\tau^{\cal{F}}$ are in $Y$ with  $|x_\tau - x_{\tau^{\cal{F}}}^{\cal{F}}| < \eps$, and since $F$ is
  Lipschitz, we have
  $$
  |F(x_\tau) - F(x_{\tau^{\cal{F}}}^{\cal{F}})| <\eps
  \mathrm{Lip}_Y^\eps F.
  $$

  Case II. If $\tau < \tau^{\cal{F}}$, then $x_\tau \in
  Y$ while $x_\tau^{\cal{F}} \not \in Y$, but still
  $d (x_\tau,x_\tau^{\cal{F}}) < \eps$. Player I continues by pulling towards $x_\tau$, and Lemma~\ref{vbound} tells us that the conditional expected value $\mathbb{E}\big[F(x_{\tau^\cal{F}})\mid x_\tau \big]$ to player I of the II-favored $\eps$-game thus played is at least
  $$
  F(x_\tau)-(2\eps + d^\eps (x_\tau,x_\tau^{\cal{F}}))
  \mathrm{Lip}_Y F \geq F(x_\tau)- 3\eps \mathrm{Lip}_Y F\,.
  $$

  Case III. If $\tau^{\cal{F}} < \tau$, then $x_{\tau^{\cal{F}}}^{\cal{F}} \in
  Y$ and $x_{\tau^{\cal{F}}} \not \in Y$. Let $y=x_{\tau^{\cal{F}}}^{\cal{F}}$ and
  $z=x_{\tau^{\cal{F}}}$; recall that $d(y,z) < \eps$.  By Lemma \ref{ubound}, we have
  $$
  u_I^\eps(z) \leq F(y) + O(1)\, e^{\b\,\mathrm{diam} X} \, \eps \, \mathrm{Lip}_Y F.
  $$
 Thus $u_I^\eps(z) \leq F(y) + O(\eps)$. Continuing from $z$, player I cannot do better than $u_I^\eps(z)$, hence $\mathbb{E}\big[F(x_\tau) \mid x_{\tau^\cal{F}}\big] \leq
  F(x_{\tau^\cal{F}}) + O(\eps)$.

Summarizing the three cases, $\tau^\cal{F}<\infty$ a.s., and
$$\mathbb{E} F(x_{\tau^\cal{F}}) + O(\eps) \geq \mathbb{E} F(x_{\tau}),$$
hence
$ F_-(\S_I^\cal{F},\S_{II}^\cal{F}) \geq F_-(\S_I,\S_{II}) - O(\eps)$, as desired.
\end{proof}

\begin{proposition}
  $\|w^\eps - u^\eps\|_\infty = O(\sqrt{\eps})$.
\label{wconv}
\end{proposition}

\begin{proof}
  The proof is identical to the proof of Proposition~\ref{vconv}, except that for case II
  we use $w^\eps \leq F(y) + O(\sqrt{\eps})$ from Lemma~\ref{wbound}, and for case III we
  use $u^\eps \geq v^\eps \geq F(y) -O(\eps)$ from Lemma~\ref{vbound}.
\end{proof}

\begin{corollary}
If the odds function $\rho(\cdot)$ satisfies (\ref{e.rhocond}), then, for any $\eps>0$, the following limits exist pointwise and are equal:
$$\vlim_{,\rho} := \lim_{n \to \infty} v^{\eps/{2^n}} = \lim_{n \to \infty} u^{\eps/{2^n}} = \lim_{n \to \infty} w^{\eps/{2^n}}\,.$$
\label{pwiselim}
\end{corollary}

\begin{proof}
By Lemma~\ref{monot}, either $(v^{\eps/2^n})_n$ or $(w^{\eps/2^n})_n$ is monotone. On the other hand, by Lemmas~\ref{vbound} and~\ref{wbound}, both sequences are bounded. Therefore, at least one of these two limits exist pointwise. On the other hand, by Propositions~\ref{vconv} and~\ref{wconv}, we have $\| v^{\eps/2^n} -w^{\eps/2^n} \|_\infty \rightarrow 0$ as $n \to \infty$. Therefore, both sequences must converge to the same pointwise limit. This is also the limit for $u^{\eps/2^n}$, since it is sandwiched between those two sequences, by Lemma~\ref{monot}.
\end{proof}

Note that the uniform closeness of $v^{\eps/2^n}$ and $w^{\eps/2^n}$ does not imply automatically that the limits in Corollary~\ref{pwiselim} hold in the sup-norm. Furthermore, at this point we do not know yet if for different $\eps$ values and $\rho(\cdot)$ functions the limits $\vlim_{,\rho}$ are the same.

%
%
\section{Comparison with exponential cones}\label{s.comp}

Recall Definitions~\ref{defEC} and~\ref{defCEC} from the Introduction. We will now show that the value function $u^\eps$ almost satisfies comparison with $\db$-exponential cones, and then that any uniform limit satisfies CEC. We may assume $\b>0$. We will write the exponential cones~(\ref{expcone}) centered at $x_0$ in the form $\varphi(x)=C(d(x,x_0))$, with $C(r)=ae^{-\db r}+b$ (positive cones) or $C(r)=ae^{\db r}+b$ (negative cones), where $a\leq 0$ and $b\in \mathbb{R}$.

\begin{lemma}
\label{almostCECA}
Let $\theta(\cdot)$ satisfy $\theta(\eps)\leq \b\eps/2 + O(\eps^2)$. Let $z\in V\subset\subset X\setminus Y$ and write $V_\eps=\{ x:\overline{B_\eps(x)} \subset V \} $ for $\eps>0$ small enough.
Let $\varphi(x) = C (d(x,z))$ be a $\db$-exponential cone centered at $z$.
If $u^\eps \leq \varphi$ on $(V \setminus V_\eps) \cup \{z\}$,
then $u^\eps \leq \varphi+O(\eps M/s)$ on $V_\eps$, where $M:=\sup_Y F-\inf_Y F$ and $s:=d(z,\p V)$.

\end{lemma}

\begin{proof}
  Starting the game at $x_0\in V_\eps$, we need to give a strategy $\S_{II}$ for player II such that for any strategy $\S_{I}$ of player I we have $F_+(\S_I,\S_{II}) \leq \varphi (x_0)+O(\eps M/s)$.

  First we suppose that $\varphi$ is a positive cone. Note that $\inf_Y F\leq u^\eps(z) \leq \varphi(z)$, and that we may assume that  $\varphi(z) \leq \sup_Y F$, since otherwise the result is trivial. We can replace $\varphi$ by a new cone that has the same value at $z$ but is otherwise as small as possible (while maintaining $u\leq \varphi$ on $(V\setminus V_\eps) \cup \{z\})$; proving the result for this $\varphi$ will suffice. Clearly, the value of $\varphi$ on $\p B_{s-\eps}(z)$ is at most $\sup_Y F$. Now, writing $C(r)=ae^{-\db r}+b$ for this $\varphi$, since the difference between the values at $z$ and $\p B_{s-\eps}(z)$ is at most $M$, we have $|a|(1-e^{-\beta (s-\eps)}) \leq M$, hence $|C'(0)|= |a|\b = O(M/s)$, after choosing $\eps<s/2$ and using $s\leq \mathrm{diam}\, X$. (Recall that $O(\cdot)$ is allowed to depend on $X$.)

Fix $\delta >0$, and let $x_{k-1} \in V_\eps$ and $r=d(x_{k-1},z)$. Let $\S_{II}$ be the
  pull toward $z$ strategy: if
  $r\geq \eps$, then player II reduces the distance to $z$ by
  almost $\eps$, enough so that $C(d(x_k,z)) < C(r-\eps) +
  \delta/2^k$; if $r<\eps$, then player II moves to $z$; note that
$C(0) < C(r-\eps)+O(|a|\b\eps)=C(r-\eps)+O(\eps M/s)$ in this case. Then,
\begin{eqnarray*}\mathbb{E}[\varphi(x_k) | {\cal{F}}_{k-1}] &<& \cefrac \,
  C(r+\eps) + \fracce \left( C(r -\eps) + \frac{\delta}{2^k} + O(\eps M/s)\,\1_{\{r<\eps\}} \right) \\ &<&
  a e^{-\db r} \left(\cefrac e^{-\db\eps} + \fracce e^{\db\eps}\right) + b +
  \frac{\delta}{2^{k+1}} + O(\eps M/s)\,\1_{\{r<\eps\}}\,.
\end{eqnarray*}
By our assumption on $\theta(\eps)$, we have $\cefrac e^{-\db\eps} +
  \fracce e^{\db\eps} \geq 1-O(\eps^3)$. Recalling $a\leq 0$, this gives us
\begin{eqnarray}
  \mathbb{E}[\varphi(x_k) | {\cal{F}}_{k-1}] &<&  a e^{-\db r} + b +
  \frac{\delta}{2^{k+1}}  - a e^{-\db r}\,O(\eps^3) + O(\eps M/s)\,\1_{\{r<\eps\}}\nonumber\\
&<& \varphi(x_{k-1}) + \frac{\delta}{2^{k+1}} + O(\eps^4 M/s) + O(\eps M/s)\,\1_{\{r<\eps\}}\,.\label{phisuperMG}
\end{eqnarray}

We let $\tau= \inf\big\{ k: x_k \not \in V_\eps\setminus\{z\} \big\}$, and, for $k\leq\tau$, define
$$
M_k:=\varphi(x_k)+\frac{\delta}{2^k}-O(\eps^3 M/s)\, k - O(\eps M/s)\,\1_{\{r<\eps\}}\big|\{j: 0\leq j<k\text{ and }d(x_j,z)<\eps\}\big|\,,
$$
where the terms $O(\eps^3 M/s)$ and $O(\eps M/s)$ are chosen larger than the similar terms in~(\ref{phisuperMG}). Hence $M_k$ is a supermartingale. In order for this $M_k$ to be useful, we need that we do not lose much with these error terms until we reach $\tau$.  One might worry that having the weaker player pull towards a target raises some complications here, but fortunately, this is easily handled by $X$ being compact in the path metric, as follows.

Due to the strategy $\S_{II}$, the number of $j$'s with $d(x_j,z)<\eps$ is stochastically bounded by an exponential random variable with a uniformly bounded expectation, as long as $(1-\theta(\eps))/2$ is bounded away from $0$, which we may assume, since $\b>0$ is fixed and $\eps\to 0$.
Moreover, $\mathbb{E} [\tau] = O(\eps^{-2})$ under this strategy, by the following argument. The bias is only $\theta(\eps)=O(\eps)$, each step of player I increases the distance to $z$ by at most one step, and the distance from any $x\in X\setminus Y$ to $z$ is at most $O(1/\eps)$ steps, therefore the time $\tau$ is stochastically dominated by the time it takes for a random walk $\theta(\eps)$-biased to the right on the line segment $\{0,1,\dots,O(1/\eps)\}$ to get from some vertex $i$ to 0. This has expectation $O(\eps^{-2})$, uniformly in the starting point $i$, as follows e.g.~from the commute time formula
 of electrical networks, see \cite[Section 2.6]{LPbook} or \cite[Section 4]{Lovasz}. In our setting, if we consider the network on $\{0,1,\dots,O(1/\eps)\}$ with conductances $C(j,j+1)=\rho^j$, then the effective resistance between vertex $0$ and $i$ is  $R(0\leftrightarrow i)=\sum_{j=0}^{i-1} \rho^{-j}$, and then the formula says  
$$
\mathbb{E}_i \tau_0 + \mathbb{E}_0 \tau_i = 2 R(0\leftrightarrow i) \sum_{j=0}^{O(1/\eps)} C(j,j+1) 
= 2\, \sum_{j=0}^{i-1} O(1) \, \sum_{j=0}^{O(1/\eps)} O(1) = O(\eps^{-2})\,.
$$ 
Since we have this uniformly in $i$, we also get $\mathbb{P}[\tau > k \eps^{-2}] < \exp(-ck)$ for some $c>0$.

On the other hand, $\varphi(x_k)$ is bounded. Thus $M_k$ is uniformly integrable, and an Optional Stopping Theorem \cite[Section 4.7]{Durrett} gives
$$
 \mathbb{E}[\varphi(x_\tau)] < \mathbb{E}[M_\tau] + O\left((\eps^{3-2} + \eps) M/s \right) \leq \mathbb{E}
 [M_0] + O(\eps M/s) = \varphi(x_0) +  \delta + O(\eps M/s)\,.
$$
Since $\tau<\infty$ a.s.~is clear from $\S_{II}$, we have $x_\tau\in (V\setminus V_\eps) \cup \{z\}$, so $u^\eps(x_\tau)\leq \varphi(x_\tau)$, hence $u_{I}^\eps (x_0) \leq \sup_{\S_{I}} F_-(\S_I, \S_{II}) \leq \sup_{\S_{I}} \mathbb{E} u_I^\eps(x_\tau) \leq  \sup_{\S_{I}} \mathbb{E} \varphi (x_\tau) \leq \varphi (x_0) + \delta + O(\eps M/s)$ for any $x_0\in V_\eps$. Since $\delta$ was arbitrary, we have $u^\eps = u_{II}^\eps \leq  \varphi + O(\eps M/s)$ on $V_\eps$, as claimed.

If $C(r)=ae^{\db r}+b$, then player II pulls away from $z$ and the same equations result.
\end{proof}

\begin{lemma}
\label{almostCECB}
Let $\theta(\cdot)$ satisfy $\theta(\eps) \geq \b\eps/2 - O(\eps^2)$. Let $z\in V\subset\subset X\setminus Y$ and
  write $V_\eps=\{ x:\overline{B_\eps(x)} \subset V \} $ for $\eps>0$ small enough. Let $\varphi(x):=C(d(x,z))$ be a $\db$-exponential cone centered at $z$. If $u^\eps \geq \varphi$ on $(V \setminus V_\eps) \cup \{z\}$, then $u^\eps \geq \varphi-O(\eps M/s)$ on $V_\eps$, where $M:=\sup_Y F - \inf_Y F$ and $s:=d(z,\p V)$.
\end{lemma}

\begin{proof} It is the same as the proof of Lemma~\ref{almostCECA}, we just have to reverse the roles of players I and II and the direction of all the inequality signs. The bound $\mathbb{E}[\tau]=O(\eps^{-2})$ is even easier now, since the stronger player is pulling towards a target.
\end{proof}

We want to use the above approximate CEC lemmas for $u^\eps$ to prove that if they converge to some $u$, then $u$ satisfies CEC exactly. But in order to be able to apply these lemmas,
we need that if $u^\eps \leq \varphi$ (or $u^\eps \geq \varphi$) on $\partial V$ for $V\subset\subset X\setminus Y$, then this almost holds also in $V\setminus V_\eps$. This continuity property is actually a simple consequence of the lemmas themselves:

\begin{lemma}[Uniform Lipschitz away from $Y$] Assume that $|\theta(\eps)-\b\eps/2|=O(\eps^2)$. Let $x_1,x_2\in X\setminus Y$ with $\max\{d(x_i,Y)\}=s$. Recall $M=\sup_Y F - \inf_Y F$. Then, for any $\eps>0$ with $2\eps<\min\{d(x_i,Y)\}$, we have $|u^\eps(x_1)-u^\eps(x_2)| \leq O\left(d^\eps(x_1,x_2)\, M/s\right)$.
\label{uniflip}
\end{lemma}

\begin{proof}
Assume $d(x_1,Y)=s$. Let $V=(X\setminus Y)_\eps$, and take the positive $\db$-exponential cone $\varphi(x)=C(d(x_1,x))$ centered at $x_1$ with $\varphi(x_1)=u^\eps(x_1)$ and $\inf\{\varphi(y): y\in X\setminus V_\eps\}=M$. Then $C'(0)=O(M/s)$, while, by Lemma~\ref{almostCECA}, we have
\begin{eqnarray*}
u^\eps(x_2) \; \leq \; \varphi(x_2)+O(\eps M/s) &\leq& \varphi(x_1)+C'(0)\,d(x_1,x_2)+O(\eps M/s)\\
&\leq& u^\eps(x_1)+O\left(d^\eps(x_1,x_2)\, M/s\right).
\end{eqnarray*}

Similarly, from looking at a negative cone centered at $x_1$ and Lemma~\ref{almostCECB}, we obtain that $u^\eps(x_2)\geq u^\eps(x_1)-O(d^\eps(x_1,x_2)\, M/s)$, and we are done.
\end{proof}

Notice that we can now plug the Lipschitz bound of Lemma~\ref{uniflip} into the proofs of Lemmas~\ref{almostCECA} and~\ref{almostCECB} to get the same results with $V\subset\subset X\setminus Y$, any $z\in X\setminus Y$ inside or outside $V$, and $s=d(z,Y)$.
\medskip

With the help of the regularity provided by Lemma~\ref{uniflip}, we can now show that the convergence in Corollary~\ref{pwiselim} is in fact uniform:

\begin{theorem}
If the odds function $\rho(\cdot)$ satisfies (\ref{e.rhocond}), then, for any $\eps>0$, we have the uniform limits
$$\vlim_{,\rho}= \lim_{n \to \infty} v^{\eps/{2^n}} = \lim_{n \to \infty} u^{\eps/{2^n}} = \lim_{n \to \infty} w^{\eps/{2^n}}\,.$$
\label{uniflim}
\end{theorem}

\begin{proof}
Using Lemmas~\ref{vbound} and~\ref{wbound} at $Y$, and Lemma~\ref{uniflip} inside $X\setminus Y$, we get that the sequence $(u^{\eps/2^n})_n$ is uniformly bounded and ``asymptotically uniformly equicontinuous'': $\forall\,\eta>0$ $\exists\,\delta>0$ and $N\in\mathbb{N}$ such that whenever $n>N$ and $d(x,y)<\delta$ for some $x,y\in X$, we have $\big|u^{\eps/2^n}(x)-u^{\eps/2^n}(y)\big| < \eta$. The proof of the Arzel\`a-Ascoli theorem goes through for such sequences, and our underlying metric space $X$ is compact, hence there is a uniformly convergent subsequence of $(u^{\eps/2^n})_n$. Since the sequence has a pointwise limit, this limit must be uniform, as well. By Propositions~\ref{vconv} and~\ref{wconv}, the limits for $(v^{\eps/2^n})_n$ and $(w^{\eps/2^n})_n$ are also uniform.
\end{proof}

We can also derive now the key property of subsequential uniform limits along any subsequence $\eps_n\searrow 0$:

\begin{theorem}
Any subsequential uniform limit $\tilde{u}= \lim_{n \to \infty} u^{\eps_n}$ satisfies comparison with $\db$-exponential cones in $X$, and hence it is continuous on $X$.
\label{compcones}
\end{theorem}

\begin{proof}
Consider an open set $V \subset\subset X\setminus Y$ with $d(V,Y)=s$, and a $\db$-exponential
  cone $\varphi$ centered at $z\in X$ such that $\tilde{u} \leq \varphi$ on $\partial (V\setminus \{z\})$. We must show that $\tilde{u} \leq \varphi$ on $V$.
Fix $\delta >0$. For $n$ large enough, $u^{\eps_n} \leq \tilde{u} +
  \delta \leq \varphi + \delta$ on $\partial V$. For all $ x \in V
  \setminus V_{\eps_n}$, there exists a $ y \in \partial V$ such
  that $|x-y| < \eps_n$. Note that $\varphi$ is uniformly continuous on $\overline{V}$,
hence $\varphi(y)\leq\varphi(x)+\delta$ if $n$ is large enough. Also, by
Lemma~\ref{uniflip}, we have $u^{\eps_n}(x) \leq u^{\eps_n}(y) + O(\eps_n \cdot M/s)$, which is at most $u^{\eps_n}(y)+\delta$ for sufficiently large $n$. Summarizing, we have $u^{\eps_n}(x) \leq u^{\eps_n}(y) +\delta \leq \varphi(y) + 2\delta \leq \varphi(x)+3\delta$ for $x\in V \setminus V_{\eps_n}$.

By the extension of Lemma \ref{almostCECA} mentioned after Lemma~\ref{uniflip}, $u^{\eps_n}$ almost satisfies comparison with the $\db$-exponential cone $\varphi+3\delta$ from above, hence $u^{\eps_n} \leq \varphi + 3 \delta + O(\eps_n \cdot M/s)$ on $V$. Since $\delta$ was arbitrary, we have $\tilde{u} = \lim_{n\to \infty} u^{\eps_n} \leq \varphi$ on $V$.

The proof that $\tilde{u}$ satisfies comparison with $\db$-exponential cones from below is the same, except that the reference to Lemma~\ref{almostCECA} has to be replaced by Lemma~\ref{almostCECB}.

That CEC implies the continuity of $\tilde{u}$ at points $x\in X\setminus Y$ is clear by following the proof of Lemma~\ref{uniflip}. Continuity at points $x\in Y$ follows from Lemmas~\ref{vbound} and~\ref{wbound}.
\end{proof}

%
%
\section{Uniqueness}\label{s.unique}

Let us give an outline of the strategy for proving uniqueness, which will be somewhat similar to \cite[Section 3.5]{PSSW2}. If $u$ is a continuous extension of $F$ that satisfies comparison with $\db$-exponential cones from above, then, provided that $\theta(\eps)$ is large enough, player I can more-or-less achieve a payoff $u$ in the $\eps$-game. Similarly, if $u$ satisfies CECB and $\theta(\eps)$ is small enough, then player II can more-or-less achieve $u$ in the $\eps$-game. Putting these two parts together, we get that if $u$ satisfies CEC, and the bias is exactly $\theta_0(\eps)=\tanh(\b\eps/2)$, then both players can almost achieve $u$; more precisely, we get $v^\eps < u < w^\eps$. This implies that $u=\vlim_{,\rho_0}$, for any $\eps$. Therefore, all continuous CEC extensions of $F$ are the same.

\begin{lemma}
  Let $w : X \longrightarrow \R$ be continuous and satisfy comparison with $\db$-exponential cones from above on $X \backslash Y$. Suppose that $\theta(\eps)$ satisfies $\frac{1}{1+e^{-\db\eps}}\leq \frac{1+\theta(\eps)}{2}$, and fix $\delta >0$. Then, in I-favored tug-of-war, player I may play to make $M_{k \wedge \tau}$ a submartingale, where $M_k:= w(x_k) -\delta/{2^k}$ and $\tau := \inf \{k: d(x_k,Y) < 3\eps \}$.
\label{unique1}
\end{lemma}

\begin{proof}
  Let $X^{3\eps}:= \{ x \in X : d(x,Y) \geq 3\eps \} $, and so $\tau$ is the time to exit $X^{3\eps}$. Let $k\leq \tau$; then $x_{k-1} \in X^{3\eps}$, and for player II's proposal $z=z_k$ we have $B_{2\eps}(z) \subset X \backslash Y$.

We need to show that  $\mathbb{E}[M(x_k)| z_k ,
  {\cal{F}}_{k-1}] \geq M(x_{k-1})$. We are not assuming that suprema and infima of $w$ over balls are achieved; however, player I may play so that $\mathbb{E}[w(x_k)| z_k ,
  {\cal{F}}_{k-1}] > \cefrac \sup_{y\in B_{2\eps}(z)}w(y) +\fracce w(z) - \frac{\delta}{2^k}$.
We let
$$
p:=\sup_{y\in B_{2\eps}(z)}w(y) \qquad\text{and}\qquad q :=\cefrac\, p +\fracce\, w(z)\,.
$$
We will show using CECA that for all $x\in B_\eps(z)$ we have $q \geq w(x)$. Since $x_{k-1} \in B_\eps(z)$, this will imply
that $\mathbb{E}[w(x_k) | z_k , {\cal{F}}_{k-1}] >  q  -  \frac{\delta}{2^k} \geq w(x_{k-1}) - \frac{\delta}{2^{k}}$ and $\mathbb{E}[M(x_k)| z_k ,
  {\cal{F}}_{k-1}] > w(x_{k-1}) - \frac{\delta}{2^k} - \frac{\delta}{2^k} = M_{k-1}$, as desired.

We know that $w(z) \leq \sup_{y\in B_{2\eps}(z)}w(y)=p$. First we notice that if $w(z) = p$ then $q = p$ as well, and hence $w(x) \leq \sup_{y\in B_{2\eps}(z)}w(y) =q$ for all $x \in B_{\eps}(z)$, as we wished.

Next we assume that $p>w(z)$. We let $V:=B_{2\eps}(z) \backslash \{ z \}$. We will define a positive exponential cone $\varphi(x)=C(d(z,x))$ centered at $z$, such that $w(x) \leq \varphi(x)$ on $x\in\partial V$. By CECA we will have that this inequality holds on all of $V$. In addition, the exponential cone will be constructed so that $\varphi(x) \leq q$ on $x\in\partial B_{\eps}(z)$. Since $\varphi(x)$ is increasing in distance from $z$, this will imply that for $x\in  B_\eps(z)$ we have $w(x) \leq \varphi(x) \leq q$, as desired.

 To that end, we let $C(r) := \frac{p-w(z)}{1-e^{-2\b\eps}} (1-e^{-\db r}) + w(z)$. Notice that $C(0)=w(z)$ and $C(2\eps) = p$, hence $w(x) \leq C(d(x,z))$ on $x\in\partial V$. Furthermore,
$C(\eps)= \frac{p-w(z)}{1+e^{-\db\eps}} + w(z) \leq \cefrac (p-w(z)) + w(z) =q$, by our assumption on $\theta(\eps)$. Thus $\varphi$ has all the desired properties, and player I can play so that  $\mathbb{E}[M_k| z_k , {\cal{F}}_{k-1}] \geq M_{k-1}$.
\end{proof}

\begin{lemma}
  Suppose $w : X \longrightarrow \R$ is continuous, satisfies
  comparison with $\db$-exponen\-tial cones from above on $X \backslash
  Y$, and $w\leq F$ on $Y$. Let $\theta(\eps)$ satisfy $\frac{1}{1+e^{-\db\eps}}\leq \frac{1+\theta(\eps)}{2}$. Then $w^\eps \geq w$ for all $\eps>0$, where $w^\eps$ is the value of the I-favored $\eps$-game for player II.
\label{unique2}
\end{lemma}

\begin{proof}
 Fix $\delta >0$, and let $M_k$, $\tau$ and $X^{3\eps}$ be as in Lemma~\ref{unique1}. We know by that lemma that player I can play such that $\mathbb{E}[M_{k\wedge \tau}] \geq M_0$. Since we want to show that $\inf_{\S_{II}}\sup_{\S_I} F_+(\S_I,\S_{II})\geq w$, we may assume that $\S_I$ is the above strategy of player I, while $\S_{II}$ is a strategy that ensures that the game almost surely finishes (which exists by Lemma~\ref{small}). In particular, we may assume that $\tau<\infty$ a.s.

Now, $w^\eps(x_0) \geq \mathbb{E} [w^\eps(x_\tau)-\delta 2^{-\tau}]$ if $x_\tau$ is reached under the above $\S_I$ and an almost optimal $\S_{II}$. On the other hand, $w^\eps(x_{\tau}) -\delta 2^{-\tau}=w(x_{\tau}) - \delta 2^{-\tau}+\big(w^\eps(x_{\tau}) - w(x_{\tau})\big) \geq M_{\tau} - \lambda_\eps$, where $\lambda_\eps := \sup_{X \backslash X^{3\eps}} (w - w^\eps)$. Now recall that $\mathbb{E}[M_{\tau}] \geq M_0= w(x_0) - \delta$. Altogether, we get that $w^\eps(x_0) > w(x_0) - \lambda_\eps - \delta$.

Letting $\delta=\eps$, we have $w^\eps \geq w - \lambda_\eps -
  \eps$. We will show in the next paragraph that $\limsup_{\eps \searrow 0} \lambda_\eps
  \leq 0$, hence letting $\eps \to 0$ gives that $\limsup_{\eps
    \searrow 0} w^\eps \geq w$. We know that $w^\eps \leq
  w^{2\eps}$, and so $w^\eps \geq \limsup_{\eps ' \searrow
    0} w^{\eps '} \geq w$, as desired.

Given any $x\in X\setminus X^{3\eps}$, let $y \in Y$ be such that $d(x,y)\leq 3\eps$. Notice that $a:=(\mathrm{Lip}_Y F) \sup_{y,y'\in Y} \frac{d(y,y')}{1-e^{-\db d(y,y')}} < \infty$ because $Y$ is compact. Now let $C(r):=a(1-e^{-\db r})+F(y)$, and $\varphi(x):=C(d(x,y))$, a $\db$-exponential cone centered at $y$. Then, for $y' \in Y$ we have $\varphi(y') \geq (\mathrm{Lip}_Y F) d(y,y')+ F(y) \geq F(y') \geq w(y')$, where the last inequality is by assumption. This implies
  $\varphi \geq w$ on $X$ by comparison with exponential cones
  from above. Thus $ w(x) \leq \varphi (x) \leq a(1-e^{-3\b\eps}) + F(y)$. On the other hand,
$w^\eps(x) \geq F(y) - O(\eps)$, by $w^\eps(x)\geq v^\eps(x)$ and Lemma~\ref{vbound}.
Therefore, $w(x)-w^\eps(x) < a(1-e^{-3\b\eps}) +O(\eps)$ and
  $\limsup_{\eps \searrow 0} \lambda_\eps \leq 0$, as promised above.
\end{proof}

For functions satisfying CECB instead of CECA, we have the following analogue:

\begin{lemma}
  Suppose that $v : X \longrightarrow \R$ is continuous, satisfies
  comparison with $\db$-exponential cones from below on $X \backslash
  Y$, and $v\geq F$ on $Y$. Let $\theta(\eps)$ satisfy $\frac{1}{1+e^{-\db\eps}}\geq \frac{1+\theta(\eps)}{2}$. Then $v^\eps \leq v$ for all $\eps>0$, where $v^\eps$ is the value of the II-favored $\eps$-game for player I.
\label{unique3}
\end{lemma}

\begin{proof}
The analogue of Lemma~\ref{unique1} holds exactly, and then the proof of Lemma~\ref{unique2} can also be modified, the only difference being that the reference to Lemma~\ref{vbound} giving $w^\eps(x) \geq v^\eps(x) \geq F(y) - O(\eps)$ has to be replaced by Lemma~\ref{wbound} giving $v^\eps(x) \leq w^\eps(x) \leq F(y) + O(\sqrt{\eps})$.
\end{proof}

We can now put together the pieces to prove our main result:

\begin{proof}[Proof of Theorem~\ref{t.general}] By Theorem~\ref{uniflim}, if $\eps>0$ and $\rho(\cdot)$ satisfies~(\ref{e.rhocond}), then the sequences $u^{\eps/2^n}$, $v^{\eps/2^n}$, $w^{\eps/2^n}$ all converge uniformly to some $\vlim_{,\rho}$ as $n\to\infty$. By Theorem~\ref{compcones}, this $\vlim_{,\rho}$ satisfies comparison with $\db$-exponential cones, and is clearly an extension of $F$. On the other hand, if $u$ is an extension that satisfies CEC, then playing the favored $\eps$-games with $\theta_0(\eps)=\tanh(\b\eps/2)$, Lemmas~\ref{unique2} and~\ref{unique3} give that $v^\eps\leq u \leq w^\eps$ for any $\eps>0$, hence $u=\tilde{u}_{\eps,\rho_0}$, for any $\eps>0$. Consequently, all extensions $u$ (in particular, all limits $\vlim_{,\rho}$) are equal to each other. Finally, we get the uniform full convergence $\|u^\eps-u\|_\infty\to 0$ using the argument employed in Theorem~\ref{uniflim}: the entire family $(u^\eps)_{\eps\searrow 0}$ is uniformly bounded and asymptotically uniformly equicontinuous, hence any subsequence has a uniform limit; but we already know that there is a unique uniform limit, so the entire sequence converges.
\end{proof}

%
%
\section{Equivalence between viscosity solutions and CEC}\label{s.equiv}

We are going to study the {$\b$-biased infinity Laplacian equation}~(\ref{biaslap}) in $\R^n$, which may alternatively be written as
\begin{equation}
\label{cih}
\Phi u:= D^2_\nu u + \b D_\nu u =0,
\end{equation}
where $\nu:=\nabla u (x)/{\| \nabla u (x)\|}$ is the unit gradient vector; for $\nabla u(x)=0$ we define $\Phi u(x):=0$. This equation is quite degenerate, and it is not clear that the functions that one would like to call solutions are actually differentiable even once. The notion of viscosity solutions is a usual solution to this problem \cite{CIL}.

\begin{definition} Let $U\subseteq \R^n$ be a domain. An upper semicontinuous function $u$ is a viscosity solution of $\Phi u \geq 0$ in $U$, or, in other words, a viscosity supersolution of (\ref{cih}), if for every local maximum point $\tilde x\in U$ of $u-w$, where $w$ is $C^2$ in some neighbourhood of $\tilde x$, we have $\Phi w(\tilde x) \geq 0$. Similarly, a lower semicontinuous function $u$ is a viscosity solution of $\Phi u \leq 0$, or viscosity subsolution of (\ref{cih}), if for every local minimum point $\tilde x\in U$ of $u-w$ we have $\Phi w (\tilde x)\leq 0$. A function is a viscosity solution of $\Phi u=0$ if it is continuous and both a viscosity supersolution and subsolution.
\end{definition}

We now prove Theorem~\ref{t.equiv}. The idea of such a statement, together with the strategy of the proof, come from \cite{CEG}; however, some of the  details are quite different.

\begin{proof}[Proof of Theorem~\ref{t.equiv}]
We may assume that $\b>0$. We first prove that if $u$ does not satisfy CECA, then $D^2_{\nu}u+\b D_{\nu}u\geq 0$ does not hold in the viscosity sense.

Suppose we have some $V\subset\subset U$, $x_0\in U$, and a $\b$-exponential cone~(\ref{expcone}) with $\iota\in\{+,-\}$, such that $u(x)\leq C^\iota_{x_0}(x) + B \equiv C(\|x-x_0\|)$ on $x\in \partial (V\setminus \{x_0\})$, but $\exists\, \tilde x\in V$ with $u(\tilde x)>C^\iota_{x_0}(\tilde x) + B$. Then, for any smooth nonnegative function $f:U \longrightarrow \mathbb{R}$, if $\eps>0$ is small enough, then $w(x):=C^\iota_{x_0}(x)+B+\eps f(x)$ is a smooth function near $\tilde x$ that satisfies $u(\tilde x)>w(\tilde x)$, while $u(x)\leq w(x)$ on $x\in \partial (V\setminus \{x_0\})$. Hence we may choose $\tilde x$ to be a local maximum of $u-w$. We want to find an $f$ so that $D^2_{\nu}w(\tilde x)+\b D_{\nu}w(\tilde x)< 0$ holds.
We will of course search for an $f$ of the form $f(x) \equiv f_0(\|x-x_0\|)$. Then, for $\|x-x_0\|=r$,
\begin{eqnarray}\label{moveconeA}
D^2_{\nu}w(x)+\b D_{\nu}w(x)=\begin{cases}
\eps\bigl(f_0''(r)+\b f_0'(r)\bigr) & \mbox{if }C'(r)+\eps f_0'(r)>0;\\
\eps\bigl(f_0''(r)-\b f_0'(r)\bigr) & \mbox{if }C'(r)+\eps f_0'(r)<0.
\end{cases}
\end{eqnarray}

We now have to distinguish between the two cases $\iota\in\{+,-\}$. Suppose first that $\iota=+$; thus $C'(r)=A\b e^{-\b r}>0$. Let $R>0$ be so large that $V\subset B_R(x_0)$, and $f_0$ be any smooth positive function with $f_0'(r)<0$ and $f_0''(r)<0$ for all $r\in (0,R)$. If $\eps>0$ is small enough, then $C'(r)+\eps f_0'(r)>0$ for all $r\in (0,R)$. Then the LHS of (\ref{moveconeA}) is negative at $x=\tilde x\in V$, as we wished.

The case $\iota=-$ is similar. Now $C'(r)<0$, so, if we choose $f_0$ such that $f_0'(r)>0$ and $f_0''(r)<0$, then, for $\eps>0$ small enough, $C'(r)+\eps f_0'(r)<0$ for all $r\in(0,R)$, and hence the LHS of (\ref{moveconeA}) is again negative at $x=\tilde x$.

We show similarly that if $u$ does not satisfy CECB, then it is not a viscosity subsolution of (\ref{cih}). Suppose we have $u(x)\geq C^\iota_{x_0}(x)+B \equiv C(\|x-x_0\|)$ on $x\in \partial (V\setminus \{x_0\})$, but $\exists\, \tilde x\in V$ with $u(\tilde x)<C^\iota_{x_0}(\tilde x)+B$. We now want $w(x):=C^\iota_{x_0}(x)+B-\eps f(x)$ with $f\geq 0$ and $\eps>0$ small, so that $D^2_{\nu}w(\tilde x)+\b D_{\nu}w(\tilde x)> 0$ holds. Writing $f(x)\equiv f_0(r)$,
\begin{eqnarray}\label{moveconeB}
D^2_{\nu}w(x)+\b D_{\nu}w(x)=\begin{cases}
\eps\bigl(f_0''(r)-\b f_0'(r)\bigr) & \mbox{if }C'(r)-\eps f_0'(r)>0;\\
\eps\bigl(f_0''(r)+\b f_0'(r)\bigr) & \mbox{if }C'(r)-\eps f_0'(r)<0.
\end{cases}
\end{eqnarray}

One can easily see that if $f_0$ is a smooth positive function with $\sgn \,f_0'(r)=-\iota$ and $f_0''(r)>0$ for all $r\in (0,R)$, then the LHS of (\ref{moveconeB}) is positive at $x=\tilde x\in B_R(x_0)$, as we wanted.

Turning to the other direction of the equivalence claimed in Theorem~\ref{t.equiv}, let us assume that $u$ is not a viscosity supersolution at $u(0)=0$, where $0\in U$. That is, there exists a $C^2$ function $w$ such that $u(x)\leq w(x)$ near $0$, with equality at $x=0$, while $D^2_{\nu}w(0)+\b D_{\nu}w(0)<0$. Since this inequality is invariant under rotations, we may assume that $\nabla w(0)=(p,0,\dots,0)$ where $p>0$. Then the inequality for $w$ simply reads as $w_{x_1,x_1}(0)+\b p<0$. So let $q:=w_{x_1,x_1}(0)<-\b p<0$, and fix some further values $q<\tilde q<q^*<-\b p$. Writing $x=x_1e_1+y$, with $y$ being perpendicular to the unit vector $e_1$, for the Hessian matrix $\nabla^2 w(0)$ we have
$$
\langle \nabla^2 w(0)(x_1e_1+y),x_1e_1+y\rangle \leq q x_1^2+K(|x_1|\|y\|+\|y\|^2) \leq \tilde q x_1^2+\tilde K \|y\|^2,
$$
with some constants $K$ and $\tilde K$; the first inequality holds because the Hessian gives a bounded bilinear operator, while the second holds because given $\eps>0$, there is a large enough $L$ such that $xy+y^2 < \eps x^2+Ly^2$ for all $x,y>0$. From this and the Taylor expansion of $w(x)$ at $x=0$, we get
\begin{eqnarray}\label{Tay}
u(x_1,y)\leq w(x_1,y)&\leq& px_1 +\frac{1}{2}\left( \tilde q x_1^2 + \tilde K \|y\|^2\right)+O(|x_1|^3+\|y\|^3)\nonumber\\
&\leq& px_1 +\frac{q^*}{2}x_1^2 + K^*\|y\|^2 =: G(x_1,y)
\end{eqnarray}
in a neighborhood of $(x_1,y)=(0,0)$, with some constant $K^*$.

Consider the $\b$-exponential cone $\varphi(x) \equiv C^+_{z}(x)+B := A(1-e^{-\b\|x-z\|})+B$, with a choice of the centre $z$ and the slope $A$ that gives $\nabla \varphi(0)=A\b e^{-\b\|z\|}
(-z/\|z\|)=(p,0,\dots,0)$. In particular, $z=(-z_0,0,\dots,0)$ for some $z_0 >0$, and $p=A\b e^{-\b z_0}$. Then easy calculus gives $\nabla^2 \varphi(0)=\mbox{diag}\bigl(-A\b^2e^{-\b z_0},A\b e^{-\b z_0}/z_0,\dots,A\b e^{-\b z_0}/z_0\bigr)
=\mbox{diag}(-\b p,\frac{p}{z_0},\dots,\frac{p}{z_0})$, a diagonal matrix.

We made this choice for $\varphi$ because $G(x_1,y)$ has now the same gradient $(p,0,\dots,0)$ at $(x_1,y)=(0,0)$ as $\varphi$, while its Hessian is $\nabla^2 G(0)=\mbox{diag}(q^*,2K^*,\dots,2K^*)$, with $q^*<-\b p$. Note that, in defining the cone $\varphi$, we still have a freedom of choice in $z_0$ and in $B$. We take $z_0$ smaller than the distance of $0$ from $\p U$, and also so small that $p/z_0$ is much bigger than $2K^*$, and then choose $A$ to make $p=A\b e^{-\b z_0}$. Then we take $B$ such that $\varphi(0,0) < G(0,0)$, but their difference is very small even compared to $(\b p-q^*)z_0$. Then, if we take the Taylor series of $\varphi-G$ around $(0,0)$, the first derivatives are 0, but the second derivatives, both in the $x_1$ and the other directions, are positive enough compared to how little bit negative $\varphi-G$ at $(0,0)$ is, to ensure that at the boundary of the ball of radius $z_0$ around $(0,0)$  we have $\varphi-G > 0$. That is, for $V:=B_{z_0}(0)$, we have $u(x)\leq w(x)\leq G(x) < \varphi(x)$ for $x\in\partial V$, while $u(0)=w(0)=G(0) > \varphi(0)$, hence CECA fails for $u$.

The proof that if $u$ is not a viscosity subsolution, then CECB fails, is almost identical to the above. We only need to change the direction of all the inequality signs, and thus can construct a negative cone for which CECB fails.
\end{proof}

%
%
\section{Concluding remarks and open problems}\label{s.open}

\noindent{\bf Explicit biased infinity harmonic functions.} There are only few explicit infinity harmonic functions known even in the unbiased case: in two dimensions, $u(x,y)=|x|^{4/3}-|y|^{4/3}$ is one example, the argument function in a sector (in polar coordinates, $u(r,\phi)=\phi$) is another. See \cite{A3,PSSW2} for further examples. For the biased case, it remains a challenge to give any genuinely non-1-dimensional example. For instance, it is not difficult to show that in a sector, with constant 0 and 1 boundary values on the bounding half-infinite rays, no extension can be just a function of the angle $\phi$.

Even if we do not have explicit solutions, we could have an efficient algorithm to compute the unique extension numerically. In the unbiased case, using the absolutely minimizing Lipschitz property, there is a simple algorithm computing the discrete infinity harmonic extension on a finite graph, in polynomial time in the size of the graph \cite{LLPU,LLPSU}. The same idea can be used to approximate continuum solutions in $\R^n$ \cite{Ober}. However, we do not know of anything like that in the biased case.
\medskip

\noindent{\bf Uniqueness questions.} In what generality does uniqueness hold? For example, what happens if $U=X\setminus Y$ is a bounded domain in $\R^n$, but $X$ is compact only in the Euclidean metric, not in the path metric? Or is there a unique extension in a sector of any angle, with constant 0 and 1 boundary values on the bounding half-infinite rays?

When $X$ is an arbitrary length space, $Y\subset X$ is closed, and $F$ is a bounded Lipschitz function on $Y$, there might be several bounded continuous extensions satisfying CEC, but, when $\b>0$, the ``smallest'' extension seems to be more canonical than the others. Is the pointwise infimum of all bounded CEC extensions finite and also a CEC extension? (For $\b<0$ we should take the pointwise supremum.) A simple example is the half-line $X=[0,\infty)$ with $Y=\{0\}$ and $F(0)=0$, mentioned after Theorem~\ref{t.main}. There are infinitely many viscosity solutions, or equivalently, CEC extensions here, with pointwise infimum $u(x)\equiv 0$. Now, if we strengthen Definition~\ref{defCEC} for CEC to include comparison over all open subsets $V\subset \overline{V} \subset X\setminus Y$ (bounded or not), then $u(x)\equiv 0$ becomes the unique such extension. This ``best'' extension is also singled out from the point of view of the game: it is easy to see that $u_I^\eps(x)\equiv 0$ and $u_{II}^\eps(x)\equiv +\infty$ for all $\eps>0$, hence the game does not have a value in the strong sense we used in this paper, but there are some natural weaker versions according to which there is a value, $u^\eps(x) \equiv 0$. For example, in the definition of $u_{II}^\eps$, we could take the sup only over such strategies $\mathcal{S}_I$ for which there exists some $\mathcal{S}_{II}$ terminating the game almost surely, because we should not punish player II with a payoff $+\infty$ for player I forcing the game not to terminate. (There would be a similar change in the definition of $u_I^\eps$, as well.) Another solution is to assign a payoff $L\in[-\infty,\infty]$ to any non-terminating game play. This game will have a value function $u^\eps_L(x)$, which is constant $0$ for all $L\leq 0$, and some larger $L$-dependent CEC function when $L>0$, so the $L\leq 0$ case may be called canonical.

So, we have the natural questions: in what domains is it true that the pointwise infimum is the unique extension that satisfies the stronger CEC condition, moreover, is the $\eps\to 0$ limit of value functions in some or most of the weaker senses?
\medskip

\noindent{\bf Regularity of solutions.} A major open problem in the unbiased case is whether infinity harmonic functions in $\R^n$ are in the class $C^1$; see \cite{CEG,CE,EY}. For topological reasons, the understanding is better in dimension 2, where solutions are known to be in $C^{1,\alpha}$ \cite{Savin,ES}, for some $\alpha>0$. One cannot expect a much stronger result, such as $C^2$, as shown by the example $|x|^{4/3}-|y|^{4/3}$. All these proofs use comparison with cones; it would be worth examining how much goes through to the biased case, using comparison with exponential cones. Do we have better or worse regularity than in the unbiased case?

\medskip

\noindent{\bf Expected duration of game.} Given $X$ compact in the path metric, in the unbiased $\eps$-game, if both players play optimally, then the ``local variation'' $\delta^+(x_n)=\sup_{y\in B(x_n,\eps)} u^\eps(y)-u^\eps(x_n)$ is non-decreasing in $n$, and thus it is not hard to see that the expected duration of the game is $O_{X}(\eps^{-2})$. However, 
it is much less clear that a player can achieve an (almost) optimal payoff within a similar expected time, regardless of what the other player does. This is closely related to a PDE stability question; let us formulate the exact connection in the unbiased case, where more information is available.

In \cite{PSSW2}, the $\eps$-game with running payoff was introduced: in addition to our usual data $X,Y,F$, we have a continuous $f:X\longrightarrow \R$, and the payoff for player I, if the game ends at step $\tau$, is $F(x_\tau)+ \eps^2 \sum_{i=0}^{\tau-1} f(x_i)$. It was shown that if $X$ is compact, $F$ is Lipschitz on $Y$, and $\min_X f>0$ or $\max_X f<0$, then the $\eps$-game has a value $u^\eps$, which converges uniformly to a continuous function $u$ that is the unique viscosity solution of the inhomogeneous equation $\Delta_\infty u = f$, with boundary condition $u\big|_{Y}=F$.

\begin{proposition}
\label{pr.duration}
Assume that $X,Y,F$ are as in the previous paragraph. Then the following are equivalent:
\begin{enumerate}
\item For any $\delta,\eps>0$, in the standard $\eps$-tug-of-war, player I has a strategy that achieves a payoff at least $u^\eps-\delta$ starting from any position $x_0\in X\setminus Y$, and the expected duration of the game is at most $C_\delta \eps^{-2}$.
\item If $f_k$ is continuous with $\min_X f_k > 0$ for all $k\in\N$, converging uniformly to 0 as $k\to\infty$, then the viscosity solutions $u_k$ to $\Delta_\infty u_k = f_k$ converge uniformly to the solution of $\Delta_\infty u=0$, with boundary value $F$.
\end{enumerate}
\end{proposition}

\begin{proof} Assuming (1), fix a small $\delta>0$, then choose $\eta>0$ smaller than $\delta / C_\delta$. If $k$ is large enough, then $\eta > f_k > 0$. The reflected version of (1) tells us that, in the standard $\eps$-game, player II can make sure that the payoff is at most $u^\eps+\delta$, and the game finishes within expected time $C_\delta \eps^{-2}$. The exact same strategy in the game with running payoff $f_k$ guarantees that the expected total payoff is at most $u^\eps+\delta+C_\delta\eta < u^\eps + 2\delta$, independently of the starting position. On the other hand, the running payoff favors player I, hence the total expected payoff is certainly at least $u^\eps$. This, together with the result we quoted from \cite{PSSW2}, establishes (2).

Assume now (2). Let $f_k \equiv 1/k$, fix any $\delta>0$, and let $k$ be large enough such that the viscosity solution $u_k$ satisfies  $\|u_k-u\|_\infty < \delta$. Now, if $\eps>0$ is small enough, then the value function $u_k^\eps$ of the $\eps$-game with running payoff $f_k$ satisfies $\|u_k^\eps-u_k\|_\infty < \delta$, again by \cite{PSSW2}. Furthermore, if $\eps>0$ is small enough, then the value function $u^\eps$ is closer than $\delta$ to the viscosity solution $u$. Altogether, if $k$ is large and then $\eps>0$ is small, then $\|u^\eps_k-u^\eps\|_\infty < 3\delta$. Therefore, there is a strategy for player II that guarantees that the payoff in the $\eps$-game with running payoff $f_k$ is at most $u^\eps+4\delta$. Since the running payoff favors player I, the same strategy achieves a payoff at most $u^\eps+4\delta$ in the standard $\eps$-game. (Starting with the reflected version of (2), we get a good strategy for player I.) Now, what is the expected duration of this game?

Since $F$ is continuous on the compact set $Y$, and $u$ is continuous on the compact set $X$, both are bounded in absolute value, say, by $C<\infty$. Hence $|u^\eps|$ is bounded by $C+\delta$, and our strategy guarantees that the payoff in the $\eps$-game with running payoff $f_k$ is at most $C+5\delta$. But then the payoff collected from the running payoff must be at most $2C+5\delta$, hence the expected duration of this game can be at most $(2C+5\delta)\eps^{-2}k$. This implies claim (1) with a constant $C_\delta=(2C+5\delta)k$, where $k=k(\delta)$.
\end{proof}

Very recently, Lu and Wang proved part (2) with PDE methods, when $X$ is the closure of a bounded domain in $\R^n$ \cite{LuWang}. It would be interesting to see a probabilistic proof of (1). More importantly, no such bound is known for the biased case, not even $C_\delta \eps^{-\gamma}$, with any exponent $\gamma<\infty$ in place of 2.


A polynomial bound on the expected duration would be helpful, for instance, in giving a polynomial time approximation scheme to compute biased infinity harmonic functions.

Let us finally note that if we had such a bound, then we would not need to be so careful in Section~\ref{s.unique} and consider the $\eps$-game with exactly $\theta(\eps)=\tanh(\b\eps/2)$: if the game does not last long, then small errors in the super- or sub-martingality could not add up to too much. We used a similar argument in Lemmas~\ref{almostCECA} and~\ref{almostCECB}, but there we did not have to bound the time till the very end of the game.

\bigskip
\noindent{\bf Acknowledgments.}
We are indebted to Scott Sheffield for helpful conversations, especially for pointing out the importance of using $\theta_0(\eps)=\tanh(\b\eps/2)$ in Section~\ref{s.unique}.

Research of Y.~Peres and S.~Somersille was supported in part by NSF grant DMS-0605166. During this work, G. Pete was a postdoc at Microsoft Research, Redmond, and at MSRI, Berkeley, and was partially supported by the Hungarian OTKA, grant T049398.

\end{document}